\documentclass[final,1p,times,number]{elsarticle}
\usepackage[numbers]{natbib}
\usepackage[colorlinks]{hyperref}
\usepackage{ulem}
\usepackage{mathrsfs}
\usepackage{verbatim,amsfonts,amssymb, mathrsfs, amsmath,   amssymb, float,amsthm,enumitem,amsthm}
\usepackage{color}
\usepackage{epsfig}

\newtheorem{theorem}{Theorem}[section]
\newtheorem{lemma}[theorem]{Lemma}
\newtheorem{prop}[theorem]{Proposition}
\newtheorem{cor}[theorem]{Corollary}
\newtheorem{definition}[theorem]{Definition}

\newtheorem{example}[theorem]{Example}
\newtheorem{remark}[theorem]{Remark}

\usepackage{tikz}
\usetikzlibrary{arrows,shapes,positioning}
\usetikzlibrary{decorations.markings}
\tikzstyle arrowstyle=[scale=1]
\tikzstyle directed=[postaction={decorate,decoration={markings,mark=at position .65 with {\arrow[arrowstyle]{stealth}}}}]
\tikzstyle reverse directed=[postaction={decorate,decoration={markings,mark=at position .65 with {\arrowreversed[arrowstyle]{stealth};}}}]
\usepackage{subfigure}

\floatstyle{ruled}
\newfloat{algorithm}{t}{loa}
\floatname{algorithm}{Algorithm}

\def\ee{{\mathcal E}}
\def\rr{{\mathscr R}}

\def\I{{\mathbb I}}

\def\E{{\mathcal E}}
\def\H{{\mathbb H}}
\def\A{{\mathcal A}}

\def\ff{{\mathcal F}}
\def\G{{\mathcal G}}
\def\V{{\mathbb V}}
\def\Y{{\mathbb Y}}
\def\C{{\mathcal C}}

\def\I{{\mathbb I}}

\def\z{{\textbf{z}}}
\def\bu{{\textbf{u}}}
\def\sat{\hbox{\rm{sat}}}

\def\trdeg{\hbox{\rm{tr.deg}}}
\def\deg{\hbox{\rm{deg}}}

\def\ord{\hbox{\rm{ord}}}
\def\lead{\hbox{\rm{ld}}}
\def\sep{\hbox{\rm{S}}}
\def\rank{\hbox{\rm{rk}}}
\def\max{\hbox{\rm{max}}}
\def\H{\hbox{\rm{H}}}

\def\S{\hbox{\rm{S}}}
\def\I{\hbox{\rm{I}}}

\begin{document}
 \title{On the Partial Differential L\"{u}roth's Theorem }

\author{Wei Li$^{\text{\,a, b}}$ and Chen-Rui Wei$^{\text{\,a, b}}$ \smallskip\\  
$^{\text{a }}$KLMM,
   Academy of Mathematics and Systems Science,\\ 
       Chinese Academy of Sciences, Beijing 100190, China;\\ 
    $^{\text{b }}$University 
of Chinese Academy of Sciences, Beijing 100049, China\\
{liwei@mmrc.iss.ac.cn, weichenrui@amss.ac.cn}  
 }

\begin{abstract}
We study the L\"{u}roth  problem for partial differential fields. 
The main result  is  the following partial differential analog of  generalized L\"{u}roth's theorem:
Let $\ff$ be a differential field of characteristic 0 with $m$ derivation operators, 
$\textbf{u}=u_1,\ldots,u_n$ a set of differential indeterminates over $\ff$.
We prove that an intermediate differential field  $\G$ between $\ff$ and $\ff\langle \textbf{u}\rangle$ is a simple differential extension of $\ff$ if and only if 
the differential dimension polynomial of $\textbf{u}$ over $\G$ is of the form $\omega_{\textbf{u}/\G}(t)=n{t+m\choose m}-{t+m-s\choose m}$ for some $s\in\mathbb N$.
This result generalizes the classical differential L\"uroth's theorem proved by Ritt and Kolchin  in the case $m=n=1$.
We  then present an algorithm to decide whether a given finitely generated differential extension field 
of $\ff$ contained in $\ff\langle \textbf{u}\rangle$ is a simple extension, and in the affirmative case, to compute a L\"{u}roth generator. 
As an application, we solve the  proper re-parameterization problem for  unirational differential curves.
\end{abstract}

  \maketitle

\section{Introduction}

One of the most important problems in algebraic geometry is to classify all the algebraic varieties up to birational equivalence. 
Rational varieties are the simplest varieties that are birationally equivalent to affine  spaces.
Unirational varieties are the ones that can be rationally parameterized, or equivalently, there are dominant rational maps from  affine spaces to unirational varieties.  
A natural problem is to ask whether a unirational  variety is rational.
In terms of field theory, this is equivalent to ask whether all non-trivial subfields of the field of rational functions are isomorphic to the field of rational functions.
This is the question J. L\"{u}roth posted in 1861, which is  known as the L\"{u}roth problem.

The L\"{u}roth's theorem, presented by  L\"{u}roth in \cite{Luroth} in 1875, 
states that if $L$ is a subfield of a purely transcendental extension $k(t)$ of $k$, 
strictly containing $k$, then $L$ is purely tanscendental.
Geometrically, this is to say  a unirational curve over any field is always rational.
A generalization of L\"{u}roth's theorem to transcendence degree one subfields of rational function fields was proved by Igusa \cite{Igusa} and Samuel \cite{Samuel} in 1950s, and  by Ollivier and Sadik \cite{Ollivier2022}with a simple and constructive proof in 2022. 
G. Castelnuovo solved the L\"{u}roth problem positively for  rational function fields in two variables over an algebraically closed field in 1893,
and thus, a unirational surface over algebraically closed field is  rational.
However, the L\"{u}roth problem fails for unirational varieties of dimension greater than two over any field.
\vskip3pt
 
Differential algebra, or differential algebraic geometry, is a subject aiming to  generalize the theory of algebraic geometry to an analogous theory for algebraic (ordinary and partial) differential equations.
In 1932, Ritt  proved the ordinary differential version of L\"{u}roth's theorem in the case when the base  field $k$ is a differential field of meromorphic functions in an open set of the complex plane and $L$ a finitely differentially generated extension of $k$ in  \cite{Ritt1932}. 
Later, E. Kolchin  proved the general differential L\"{u}roth's theorem without finiteness restrictions in \cite{Kolchin1944, Kolchin1947}: 
\begin{quote} {\it(Differential L\"{u}roth's Theorem)} Let $k$ be an ordinary differential field of characteristic 0 and $k\langle u\rangle$ a purely differentially transcendental extension of $k$. If $L$ is a differential field with $\ff\subsetneq L\subseteq \ff\langle u\rangle$,
then there exists $v\in L$ such that $L=k\langle v\rangle$. 
Such a $v$ is called a L\"{u}roth generator of $L.$
\end{quote}
A geometric version of differential L\"{u}roth's theorem says that unirational differential varieties of dimension one (i.e., differential curves) are always rational.
In \cite[p. 164]{Kolchin1973}, Kolchin pointed out that his proof of the differential L\"uroth's theorem could be adapted to prove the generalized differential Lüroth’s Theorem where $u$ is replaced by a finite set of differential indeterminates and $L$ is of differential transcendence degree 1 over $\ff$.
Gao and Xu presented a constructive proof for the differential L\"{u}roth's theorem in \cite{Gaoxu}.
Contrary to the classical setting, Ollivier showed that the L\"{u}roth  problem fails  
for unirational ordinary differential varieties of dimension two  \cite{Ollivier}.
Effective differential L\"uroth's theorem was studied by  D'Alfonso et al. in \cite{DJS2014},
 where upper bounds for the degree and order of L\"uroth generators were obtained and applied to compute a L\"uroth generator with pure algebraic procedures.
Degree bounds for L\"uroth generators were refined by Freitag and Li  in \cite{FreitagLi}.  
 D'Alfonso et al. also studied the quantitative aspects for the generalized differential L\"uroth's theorem in \cite{ DJS2018}.

In this paper, we study the L\"{u}roth problem for partial differential fields.
As Kolchin pointed out in \cite{Kolchin1944}, the L\"{u}roth's theorem 
 ceases to hold for partial differential fields, and he gave a counter example  
 showing that for a partial differential field $\ff$ of characteristic 0 with two derivations $\delta_1, \delta_2$,
the differential subfield $\G=\ff\langle\delta_1(u), \delta_2(u)\rangle$ of the purely differentially tanscendental field $\ff\langle u\rangle$ is not a simple extension.
A natural question arises: under what conditions does the L\"{u}roth's theorem hold for partial differential field extensions?
This is the main problem we deal with in this article.
Specifically, we will give a necessary and sufficient condition for an intermediate differential field $\G$ between  $\ff$ and a purely differentially transcendental extension $\ff\langle u_1,\ldots,u_n\rangle$  to be a simple differential extension of $\ff$.
This partial differential L\"{u}roth's theorem is a generalization of both the classical and the generalized ordinary differential L\"{u}roth's theorem.
Then, based on the Wu-Ritt's zero decomposition theorem, we  give an algorithm to decide whether a given finitely generated differential extension field $\G\subset\ff\langle u_1,\ldots,u_n\rangle$ of $\ff$  
has a L\"{u}roth generator, and in the affirmative case, to compute a L\"{u}roth generator.

 The L\"{u}roth  problem is closely related to the proper re-parametrization problem. 
In \cite{Gao2003}, given a set of ordinary differential rational parametric equations $\mathcal P$, 
Gao proposed an algorithm to decide whether $\mathcal P$ is proper, and if $\mathcal P$ is not proper
and the implicit variety is of differential dimension one, to find a  proper re-parameterization of $\mathcal P$.
In \cite{Li2006}, Li generalized most   results of \cite{Gao2003} to the partial differential case,
including  providing a method to detect whether a given set of partial differential rational parametric equations $\mathcal P$ is proper.
But if $\mathcal P$ is an improper parametrization of a partial differential curve, 
the problem of proper re-parameterization is not dealt with.
We will solve this problem   by giving a method to decide whether a given improper $\mathcal P$  has a proper re-parameterization, and in the affirmative case, to compute a new parameter and a proper re-parameterization of $\mathcal P$. Geometrically, this method gives a criterion to decide whether a  unirational differential curve given by means of its parametrization is rational or not.

\vskip3pt
The paper is organized as follows. 
In Section 2,  we introduce basic notions and preliminary results. 
In Section 3,  we present the partial differential  version of L\"{u}roth's Theorem. 
In Section 4,  we treat the algorithmic aspects of the partial differential   L\"{u}roth's Theorem. 
In Section 5, we solve the proper re-parameterization problem for unirational partial differential curves.

\section{Preliminaries} 
In this section, we give some basic notation and preliminary results in differential algebra. 
For more details about differential algebra, please refer to \cite{BC, Kolchin1973, Ritt, Wu}.

Let $\ff$ be a differential field of characteristic 0 endowed with a finite set of derivation operators  
$\Delta=\{\delta_1,\ldots,\delta_m\}$,
and let $\E$ be a fixed universal differential extension field of $\ff$.
Let $\Theta$ = $\{\delta_1^{e_1}\ldots\delta_m^{e_m}\big| \,e_1, \ldots, e_m\in\mathbb{N}\}$.
 The {\it order} of $\theta=\delta_1^{e_1}\ldots\delta_m^{e_m}$ is defined as $\ord(\theta)=\sum_{i=1}^me_i$.
For ease of notation, we denote $\Theta_s:=\{\theta\in\Theta|\,\ord(\theta)=s\}$ 
and  $\Theta_{\leq s}:=\{\theta\in\Theta|\,\ord(\theta)\leq s\}$.
For an element $a\in\E$,  denote $a^{[s]}=\{\theta(a): \theta\in\Theta_{\leq s}\}$. 
Throughout the paper,  all differential fields   are {\it partial differential} fields, i.e., $m>1$,
and for simplicity, we often use the prefix ``$\Delta$-" as a synonym of ``partial differential" or ``partial differentially".

Let $\G$ be a $\Delta$-extension field of $\mathcal {F}$.
A subset $\Sigma$ of  $\mathcal {G}$  is said to be {\it $\Delta$-dependent} over $\mathcal
{F}$ if the set $(\theta\alpha)_{\theta \in \Theta,\alpha\in\Sigma}$ is algebraically dependent over $\mathcal {F}$,
and is said to be {\it $\Delta$-independent} over $\mathcal {F}$, or
 a family of {\it $\Delta$-$\ff$-indeterminates} in the contrary case.
 In the case $\Sigma=\{\alpha\}$, we say that
$\alpha$ is {\it $\Delta$-algebraic} or {\it $\Delta$-transcendental} over $\mathcal {F}$ respectively.
The {\it $\Delta$-transcendence degree} of $\mathcal G$ over $\mathcal F$, denoted by  $\Delta$-$\trdeg\,\mathcal G/\mathcal F$,
is the cardinality of any maximal subset of $\mathcal G$ which are $\Delta$-independent over $\mathcal {F}$.
And the transcendence degree of $\mathcal G$ over $\mathcal F$ is denoted by $\trdeg\,\mathcal G/\mathcal F$.

Let $\ff\{\Y\}\overset{\triangle}{=}\ff[\Theta(\Y)]$ be the $\Delta$-polynomial ring  in the  $\Delta$-$\ff$-indeterminates $\Y=\{y_1,\ldots,y_n\}$ over $\ff$. 
For a $\Delta$-polynomial $f$ in $\ff\{\Y\}\backslash\ff$, 
the order of $f$ is defined as  $\ord(f)=\min\{s\,|\, f\in\ff[Y^{[s]}]\}$.
A (prime) $\Delta$-ideal in $\ff\{\Y\}$ is a (prime) ideal which is closed under the derivation operators. 
  
By a $\Delta$-affine space $\mathbb A^n$, we mean  the set $\ee^n$. 
A {$\Delta$-variety} over $\ff$  is  $\mathbb V(\Sigma)=\{\eta\in\mathcal E^n\,|\, f(\eta)=0, \forall f\in\Sigma\}$ for some set $\Sigma\subseteq\mathcal {F}\{\Y\}$. 
 Given a $\Delta$-variety $V$ defined over $\ff$,  we denote $\mathbb{I}_{\ff}(V)$ to be the set of all $\Delta$-polynomials in $\mathcal F\{\Y\}$ that vanish at every point of $V$.  
For a prime $\Delta$-ideal $\mathcal P$, a point  $\eta\in \mathbb V(\mathcal P)$ is called a {\it generic point} of $\mathcal P$  (or  $\mathbb V(\mathcal P)$) if for any $f\in\mathcal F\{\Y\}$,  $f(\eta)=0$ implies $f\in\mathcal P$. By \cite{Kolchin1973}, a $\Delta$-ideal has a generic point if and only if it is prime.
 
\subsection{Differential  characteristic sets and differential dimension polynomial}

A  {\it ranking} on $\ff\{\Y\}$ is a total order on $\Theta(\Y)=\{\theta y_j|\,j=1,\ldots,n; \theta\in\Theta\}$ which is compatible with the derivation operators: 1) for any  $\theta y_j\in\Theta(\Y)$ and $\delta_k$, $\delta_k \theta y_j>\theta y_j$ 
and 2) $\theta_{1} y_{i} >\theta_{2} y_{j}$ $\Longrightarrow$ $\delta_k\theta_{1} y_{i} >\delta_k\theta_{2} y_{j}$
for $\theta_{1} y_{i}, \theta_{2} y_{j}\in \Theta (\Y)$.
 By convention, $1<\theta y_{j}$ for all $\theta y_{j}\in \Theta (\Y)$.

Two important kinds of rankings are often used:

    1) {\it Elimination ranking}: 
    $y_{i} > y_{j}$ $\Longrightarrow$ $\theta_1  y_{i} >\theta_2 y_{j}$ for any $\theta_1,\theta_2\in\Theta$.

    2) {\it Orderly ranking}: $k>l$  $\Longrightarrow$ for any $\theta_1\in\Theta_k$, $\theta_2\in\Theta_l$ and $i, j$,   $\theta_1  y_{i} >\theta_2  y_{j}$.

For a single $\Delta$-variable $y$, unlike the ordinary differential case where there is a unique ranking on $\ff\{y\}$,
there are uncountable number of rankings on  $\ff\{y\}$ even in the case $m=2$.
For the convenience of later description, we fix a {\it canonical ranking } $\mathscr{R}_0$ on $\ff\{y\}$: 
$\prod_{i=1}^m\delta_i^{k_i}(y) <\prod_{i=1}^m\delta_i^{l_i}(y)$ if and only if $ (\sum_{i} k_i, k_1,\ldots,k_{m})<_{lex}(\sum_{i} l_i, l_1,\ldots,l_{m})$. 
\vskip3pt

   Let $f$ be a $\Delta$-polynomial in $\mathcal {F}\{\Y\}$ and $\mathscr{R}$  a ranking endowed on it.  
  The greatest   $\theta y_j$   appearing effectively in $f$ is called the {\it leader} of $f$, denoted by $\lead({f})$.  The {\it rank} of $f$ is $\lead({f})^d$, denoted by $\rank(f)$, where  $d$ is the degree of $f$ in $\lead({f})$.  
   The coefficient of $\rank(f)$ in $f$ is called the {\it initial} of $f$, denoted by $\I_{f}$.  
   The partial derivative of $f$ w.r.t. $\lead({f})$ is called the {\it separant} of $f$, denoted by $\S_{f}$.
     For any two  $\Delta$-polynomials $f$, $g$ in $\mathcal {F}\{\Y\}\backslash \mathcal {F}$,
    $f$  is said to be of {\it lower rank} than  $g$ if either  $\lead(f)<\lead({g})$ or  $\lead(f)=\lead({g})$ and $\deg(f,\lead(f))<\deg(g,\lead(f))$.

    Let $f$ and $g$ be two $\Delta$-polynomials and $\rank(f)=(\theta y_j)^d$.  
  $g$ is said to be {\it partially reduced} w.r.t. $f$ if no proper derivatives of $\theta y_j$ appear in $g$.  
  And $g$ is said to be {\it reduced} w.r.t. $f$ if $g$ is partially reduced w.r.t. $f$ and $\deg(g,\theta y_j)<d$.  
A set of $\Delta$-polynomials $\mathcal {A}$ is said to be   {\it autoreduced} if each $\Delta$-polynomial of $\mathcal {A}$ is reduced     w.r.t.  any other element of $\mathcal{A}$.
   Every autoreduced set is finite.


Let $\mathcal {A}=A_{1},A_{2},\ldots,A_{\ell}$ be an autoreduced
   set with $\S_{i}$ and $\I_{i}$ as the separant and the initial of $A_{i}$.
For each $\Delta$-polynomial  $F$, there exists an algorithm, called {\it Ritt's algorithm of reduction} \cite{Sit}, which reduces
   $F$ w.r.t. $\mathcal {A}$ to a  $\Delta$-polynomial $R$ that is
   reduced w.r.t. $\mathcal {A}$, satisfying  
  $\prod_{i=1}^t\S_{i}^{d_i}\I_{i}^{e_i} \cdot F \equiv R\mod [\mathcal {A}],$ 
   for some $d_{i},e_{i}\in \mathbb N, i=1,2,\ldots,\ell$. 
We call $R$ the {\it remainder} of $F$ w.r.t. $\A$.
  
 Let $\mathcal J$ be a $\Delta$-ideal in $\ff\{\Y\}$.
 An autoreduced set $\mathcal {C}\subset\mathcal J$  is said to be a {\it characteristic set} of $\mathcal {J}$,
 if  $\mathcal {J}$ does not contain any nonzero element reduced w.r.t. $\mathcal {C}$. 
 Denote $\H_{\mathcal {C}}$ to be  the set of all the initials and
    separants of $\mathcal {C}$ and $\H_{\mathcal {C}}^\infty$ to be the minimal
    multiplicative set containing $\H_{\mathcal {C}}$.
    The {\it $\Delta$-saturation ideal} of $\C$ is defined to be
    $\sat(\C)\overset{\triangle}{=}[\mathcal
   {C}]:H_{\mathcal {C}}^\infty = \big\{p\in\ff\{\Y\}\big| \exists h\in H_{\mathcal
{C}}^\infty, \,\textrm{s.t.}\, hp\in[C]\big\}.$ 
The algebraic saturation ideal of $\C$ is denoted by $\text{asat}(\C)\overset{\triangle}{=}(\C):\H_{\mathcal {C}}^\infty$.
 If  $\mathcal {J}$ is prime, $\mathcal {C}$
reduces to zero only the elements of $\mathcal {J}$ and we have
$\mathcal {J}=\sat(\C)$.

An autoreduced set $\mathcal {C}$ is called {\it coherent} if whenever $A,  A'\in\mathcal C$ with $\lead(A)=\theta_1(y_j)$ and $\lead(A')=\theta_2(y_j)$ for some $y_j$, the remainder of $\S_{A'}\frac{\theta}{\theta_1}(A)-\S_{A}\frac{\theta}{\theta_2}(A')$ 
w.r.t. $\mathcal C$ is zero, where $\theta=\text{lcm}(\theta_1,\theta_2)$.\,\,(Here, if $\theta_j=\prod_{i=1}^m\delta_i^{a_{ji}}\,(j=1,2)$ and $\max(a_{1i},a_{2i})=c_i$, then $\theta=\prod_{i=1}^m\delta_i^{c_i}$ and $\frac{\theta}{\theta_j}=\prod_{i=1}^m\delta_i^{c_i-a_{ji}}$.)

An autoreduced set $\mathcal {A}=A_{1}, \ldots, A_{\ell}$ arranged in increasing ranks is called $\Delta$-irreducible if for any 
$k\,(1\leq k\leq\ell)$, there cannot exsist any relation of the form
$D_{k}A_{k}\equiv G_{k,1}G_{k,2}\mod[\A_{k-1}]$
in which $G_{k,1}$ and $G_{k,2}$ are $\Delta$-polynomials with the same leader as $A_k$ and reduced with respect to $\A_k$ while $D_k$ is of lower leader than $A_k$ and reduced with respect to $\A_{k-1}$,
where $\mathcal {A}_k:=A_{1}, \ldots,A_{k}$.
By the irreducibility theorem of Wu (see \cite[p.302]{Wu}), 
$\A$ is $\Delta$-irreducible if and only if $\A$ is a characteristic set of the prime algebraic ideal 
$\text{asat}(\A)$.  Combining with this,  Kolchin's criterion for an autoreduced set to be a characteristic set of a prime $\Delta$-ideal could be stated as follows.

\begin{prop}\cite[p.167, Lemma 2]{Kolchin1973}
If $\A$ is a characteristic set of a prime  $\Delta$-ideal $\mathcal P\subset\ff\{\Y\}$, then $\mathcal P=\sat(\mathcal A)$, $\mathcal A$ is coherent and  $\Delta$-irreducible.
Conversely, if $\mathcal A$ is a coherent  and  $\Delta$-irreducible autoreduced set of $\ff\{\Y\}$,
 then $\sat(\mathcal A)$ is a prime $\Delta$-ideal and $\mathcal A$ is a characteristic set of $\sat(\mathcal A)$.
\end{prop}

\vskip3pt
Prime  $\Delta$-ideals  whose characteristic sets consist of a single polynomial are of particular interest to us.

 \begin{lemma} \cite[Lemma 2.6]{li}\label{le-generalcomponent}
Let $\mathcal{P}$ be a prime $\Delta$-ideal in $\ff\{y_1,\ldots,y_n\}$ and  $A\in\mathcal{P}$.
Suppose $A$ constitutes a characteristic set of $\mathcal{P}$ under some   ranking $\mathscr{R}$.
Then $A$ is also a characteristic set of $\mathcal{P}$ under an arbitrary ranking.
In this case, we call $\mathcal{P}$ the {\rm general component} of $A$.
\end{lemma}

 Let  $\mathcal P$ be a prime  $\Delta$-ideal in  $\ff\{\Y\}$ with a generic point $\eta\in\mathbb A^n$.
The {\it $\Delta$-dimension of $\mathcal P$}, denoted by $\Delta$-$\dim(\mathcal P)$, is defined as  the $\Delta$-transcendence degree of $\ff\langle\eta\rangle$ over $\ff$. 
For  a characteristic set $\A$ of $\mathcal P$ under some ranking,
we use {$\lead(\A)$} to denote the set $\{\lead(F)\,|\,F\in\A\}$.
Call $y_j$ a {\it leading variable} of $\A$ if  there exists some $\theta\in\Theta$ such that 
$\theta(y_j)\in\lead(\A)$;  
 otherwise, $y_j$ is called a {\it parametric variable} of $\A$.
 The $\Delta$-dimension of $\mathcal P$ is equal to the cardinality of the set of all the parametric variables of $\A$.

 Kolchin introduced {the notion of differential dimension polynomial} for prime differential ideals,
 also known as Kolchin polynomial in literature,
which contains more quantitative information than the $\Delta$-dimension.

\begin{theorem} \cite[Theorem 2]{Kolchin1964}  \label{th-dimensionpoly}
Let $\mathcal P$ be a prime $\Delta$-ideal in $\ff\{y_1,\ldots,y_n\}$. 
There exists a numerical polynomial $\omega_{\mathcal P}(t)\in \mathbb Q[t]$ with the following properties:
\begin{enumerate}
\item[1)] For sufficiently large $t\in \mathbb N$, $\omega_{\mathcal P}(t)$ equals the dimension of $\mathcal P\cap\ff[\Y^{[t]}]$.
\item[2)] $\deg(\omega_{\mathcal P})\leq m=\text{\rm Card}(\Delta)$.
\item[3)] If we write $\omega_{\mathcal P}(t)=\sum_{i=0}^m a_i{t+i\choose i}$, then $a_m$ equals the $\Delta$-dimension of $\mathcal P$.
\end{enumerate}
\end{theorem}

The numerical polynomial $\omega_{\mathcal P}(t)$ is defined as  the {\it differential dimension polynomial} or  {\it Kolchin polynomial}  of $\mathcal P$.
Given a point $\eta\in\E^n$, the differential dimension polynomial of $\eta$ over $\ff$ is defined as  $\omega_{\eta/\ff}(t):{=}\omega_{\mathbb I_{\ff}(\eta)}(t)$.

\vskip4pt
Kolchin gave a criterion for a prime $\Delta$-ideal  to be the general component of some $\Delta$-polynomial.
\begin{lemma}  \label{lm-codim1} \cite[p.160, Proposition 4]{Kolchin1973}
Let $\mathcal{P}\subseteq\ff\{y_1,\ldots,y_n\}$ be a prime $\Delta$-ideal.
Then a necessary and sufficient condition for $\mathcal{P}$ to be the general component of some polynomial $A$ of order $s$
is that the differential dimension polynomial of $\mathcal{P}$ is equal to
$$\omega_\mathcal{P}(t)=n{t+m\choose m}-{t+m-s\choose m}.$$
\end{lemma}

{We will need the strong form of Wu-Ritt's zero decomposition theorem  (\cite{Wu}).}
\begin{theorem} \label{th-zerodecom}
There is an algorithmic procedure which permits to detect whether $\mathbb V(\Sigma/\mathcal D)=\emptyset$ for any two finite subsets $\Sigma,\, \mathcal D\subset K\{Y\}\setminus \{0\}$, and in the non-empty case,
to decompose $\mathbb V(\Sigma)$ in the following form
\begin{equation}
\mathbb V(\Sigma\big/\mathcal D)=\bigcup_{i=1}^k(\mathcal A_{i}\big/\H_{\A_i}\mathcal D)
\end{equation}
and 
\begin{equation}\mathbb V(\Sigma\big/\mathcal D)=\bigcup_{i=1}^k(\sat(\mathcal A_{i})\big/\mathcal D) \end{equation}
\noindent 
in which each $\mathcal A_i$ is a coherent and  $\Delta$-irreducible  autoreduced set.
 
\end{theorem}


\subsection{Simple differentially transcendental extensions}

Take an element $u\in\mathcal E$ which is $\Delta$-transcendental over $\ff$
and we consider the simple $\Delta$-transcendental extension  $\ff\langle u\rangle$ of $\ff.$
The elements of $\ff\langle u\rangle$ are quotients
  $$\eta=\frac{P(u)}{Q(u)}$$
with $P(u), Q(u)\in\ff\{u\}$, which can always be assumed to be in reduced form (that is, $P(u)$ and $Q(u)$ are relatively prime).
The highest of the orders of $P(u)$ and $Q(u)$ is called the order of $\eta$,
denoted by $\ord(\eta)$.

\medskip
The following result is a  partial differential analog of \cite[Lemma 6]{FL} in the ordinary differential case,
which will be used to give a necessary condition for the existence of a L\"{u}roth generator in section 3. 

\begin{lemma} \label{lm-ukolchinpol}
Let $P(u), Q(u)\in\ff\{u\}$ with $\gcd(P(u), Q(u))=1$ and $\ord(\frac{P(u)}{Q(u)})=s$.
Then $$ \omega_{u\big/\ff\langle\frac{P(u)}{Q(u)}\rangle}(t) = \binom{t+m}{m}-\binom{t+m-s}{m}$$
\end{lemma}
\proof  
Select an orderly ranking $\mathscr{R}$ on $\Theta(u)$, and let  $\theta(u):=\max\{\lead(P(u)), \lead(Q(u))\}$.
 We claim that ($\star$) for each $\tau\in\Theta_{k}$ with $k>0$, 
$$\tau\Big(\frac{P(u)}{Q(u)}\Big)=\frac{P_\tau(u)}{Q(u)^{k+1}} 
\text{  with }  P_\tau(u)\in\ff\{u\}  \text{ and }  \text{rk}(P_\tau(u))=(\tau\theta(u),1).$$ 
We prove the claim by induction on $k=\ord(\tau)$.
In the case $k=1$, $\tau=\delta_i$ for some $i$ and
 $\delta_i\Big(\frac{P(u)}{Q(u)}\Big)=\frac{\delta_i(P)Q-P\delta_i(Q)}{Q^{2}}$.
Set $P_\tau=\delta_i(P)Q-P\delta_i(Q).$
If $\lead(P)\neq\lead(Q)$, it is clear that $\text{rk}(P_\tau)=(\delta_i\theta(u),1).$
Otherwise, $P_\tau=(\sep_PQ-\sep_QP)\cdot\delta_i\theta(u)+T$ with some $T\in\ff\{u\}$ involving derivatives strictly less than $\delta_i\theta(u)$. 
Since $\gcd(P,Q)=1,$ $\sep_PQ-\sep_QP\neq0$ and $\text{rk}(P_\tau)=(\delta_i\theta(u),1)$ follows.
Suppose the claim holds for $\tau\in\Theta_{\leq k-1}\,(k\geq2)$ and we consider the case $\tau\in\Theta_{k}$.
There exists $\delta_i$ and $\tau_1\in\Theta_{k-1}$ s.t. $\tau=\delta_i\tau_1$.
By the induction hypothesis, $\tau_1\Big(\frac{P(u)}{Q(u)}\Big)=\frac{P_{\tau_1}}{Q(u)^{k}} 
\text{  with }  P_{\tau_1}\in\ff\{u\}  \text{ and }  \text{rk}(P_{\tau_1})=(\tau_1\theta(u),1).$ 
Since $$\tau\Big(\frac{P(u)}{Q(u)}\Big)=\delta_i\Big(\frac{P_{\tau_1}}{Q^{k}}\Big)=\frac{\delta_i(P_{\tau_1})Q-kP_{\tau_1}\delta_i(Q)}{Q^{k+1}},$$
 $P_\tau:=\delta_i(P_{\tau_1})Q-kP_{\tau_1}\delta_i(Q)$ has rank $(\tau\theta(u),1)$ and the claim is proved.

\medskip

Denote $\mathcal G=\ff\langle\frac{P(u)}{Q(u)}\rangle$.
To show $\omega_{u/\mathcal{G}}(t)=\trdeg\,\G(u^{[t]})/\G$ has the desired form,
we try to prove that for $t$ sufficiently large, $B:=u^{[t]}\backslash\Theta_{\leq t-s}(\theta(u))$ is a transcendence basis of $\G(u^{[t]})$ over $\G$. 
Note that $\ord(\theta)=s.$
By claim ($\star$), for each $\tau\in\Theta_{k}\, (1\leq k\leq t-s)$, 
 $\tau\theta(u)=\Big(Q^{k+1}\tau(\frac{P(u)}{Q(u)})-T_\tau\Big)/I_\tau\in\G\big(\sigma(u): \sigma(u)\prec_{\mathscr R}\tau\theta(u)\big)$, where $I_\tau=\text{I}_{P_\tau}$ and $T_\tau=P_\tau-I_\tau\cdot\tau\theta(u)$.
Since $\mathscr R$ is an orderly ranking,
 it is easy to show by induction that $\Theta_{\leq t-s}\theta(u)\in\G(B,\theta(u)).$
As $\theta(u)=\max\{\lead(P(u)), \lead(Q(u))\}$, $\theta(u)$ is algebraic over $\mathcal G(B)$ and thus,
$\G(u^{[t]})\subseteq\G(B,\theta(u))$ is algebraic over $\G(B)$.

\vskip3pt
Now, it suffices to show that $B$ is algebrically independent over $\G$.
Suppose the contrary. 
Then there exists $\ell\in\mathbb N$ such that $B$ is algebraically dependent over 
$\ff\big((\frac{P}{Q})^{[\ell]}\big)$. 
Since we have $\trdeg\,\ff\big(B, (\frac{P}{Q})^{[\ell]}\big)/\ff  = \trdeg\,\ff\big((\frac{P}{Q})^{[\ell]}\big)/\ff+\trdeg\,\ff\big((\frac{P}{Q})^{[\ell]}\big)(B)/\ff\big((\frac{P}{Q})^{[\ell]}\big)\leq {\ell+m\choose m}+ \binom{t+m}{m}-\binom{t-s+m}{m}-1$,
it follows that $\trdeg\,\ff(B)\big(\frac{P}{Q})^{[\ell]}\big)/\ff (B)\leq {\ell+m\choose m}-1$.
But by claim ($\star$), the Jacobian matrix of $(\frac{P}{Q})^{[\ell]}$ w.r.t. $\Theta_{\leq\ell}\big(\theta(u)\big)$ is of full rank, so  $ (\frac{P}{Q})^{[\ell]}$ should be algebraically independent over $\ff (B)$ by Jacobian criterion \cite[p. 12, 11.4]{Lefschetz}, a contradiction.
Thus,   $B$ is a transcendence basis of $\G(u^{[t]})$ over $\G$ and 
$ \omega_{u\big/\ff\langle\frac{P(u)}{Q(u)}\rangle}(t) = \binom{t+m}{m}-\binom{t+m-s}{m}$ follows. 
\qed

\medskip

Lemma \ref{lm-ukolchinpol} implies that $u$ is algebraic over $\ff\langle P(u)/Q(u)\rangle$ if and only if $\ord(P(u)/Q(u))=0.$

\vskip5pt
In section 3, we shall need the following Kolchin's result on  extension of $\Delta$-fields.
An isomorphism $\varphi$ from a $\Delta$-field $\G$ to a $\Delta$-field $\mathcal {G}'$ is a
{\it $\Delta$-isomorphism} if $\varphi\circ \delta_i=\delta_i\circ\varphi$ for each $i$. 
By a  {\it $\Delta$-isomorphism of $\G$ with respect to $\ff$}, 
we mean a  $\Delta$-isomorphim $\phi$ from $\G$ onto a $\Delta$-field $\G^{'}$ such that
\vskip2pt
 (1)  $\G^{'}$ is a $\Delta$-field extension of $\ff$; 
\vskip2pt
 (2) $\phi$ leaves every element of $\ff$ invariant; 
\vskip2pt
 (3)  $\G$ and $\G^{'}$ have a common $\Delta$-field extension.

\begin{lemma} \cite[p. 726]{Kolchin1942}\label{Kolchin-isomorphism}
Let $\G$ be a $\Delta$-field extension of $\ff$ and $\gamma\in\G$. 
A necessary and sufficient condition that $\gamma\in\ff$ is that every $\Delta$-isomorphism of $\G$
with respect to $\ff$ leaves $\gamma$ invariant. 
A necessary and sufficient condition that $\gamma$ be a primitive element, that is,
$\G=\ff\langle\gamma\rangle$, is that no $\Delta$-isomorphism of $\G$ with respect to $\ff$ other than the identity leaves $\gamma$ invariant.
\end{lemma}

We conclude this section by giving the following result that coprimality property among polynomials does not depend on base field extensions.
It should be a known result in algebra, but as we did not find a precise reference,  
we provide a simple proof  here via using generic points. 
\begin{lemma} \label{lm-coprime}
Let $f_1,\ldots,f_\ell\in\ff[X]=\ff[x_1,\ldots,x_r]$ be  relatively prime over $\ff$.
Then  $f_1,\ldots,f_\ell$ are relatively prime over  any  field extension $\G$ of $\ff$.
\end{lemma}
\proof 
We first consider the case that $\G$ is a purely transcendental extension of $\ff$,
i.e., $\G=\ff(Y)$ for a set of algebraic indeterminates $Y$ over $\ff$.
Suppose the contrary and let  $f\in\ff(Y)[X]\big\backslash\ff(Y)$ be a common divisor of  $f_1,\ldots,f_\ell$ over $\ff(Y)$.
Then there exist  $g_i\in\ff(Y)[X]$ such that $f_i=f\cdot g_i$ for each $i$.
Suppose $q(Y), q_i(Y)\in\ff[Y]$ are polynomials of minimal degree such that $qf, q_ig_i\in\ff[Y,X]$.
Then each $qq_if_i=(qf)\cdot (qM_ig_i)$ is a factorization in $\ff[Y][X].$
By Gauss's Lemma, $qq_i\in\ff$. 
Since $f_i\in\ff[X]$,   $f, g_i\in\ff[X]$ follows,  a contradiction.
Thus, $f_1,\ldots,f_\ell$ are  relatively prime over $\ff(Y)$.

For the general case, let $Y$ be a transcendence basis of $\G$ over $\ff$ and set $\G_0:=\ff(Y)$.
By the above arguments, $f_1,\ldots,f_\ell$ are  relatively prime over $\G_0$.
Suppose the contrary and let  $g\in\G[X]\big\backslash\G$ be a common irreducible divisor of  
$f_1,\ldots,f_\ell$ over $\G$. 
Let $\eta$ be a generic point of $(g)_{\G[X]}$.
Consider the vanishing ideal $\mathcal P=\mathbb I_{\G_0}(\eta)$ of $\eta$ over $\G_0$.
Since $\G$ is algebraic over $\G_0$ and $\trdeg\, \G(\eta)/\G=n-1$,
we have $\trdeg\, \G_0(\eta)/\G_0=n-1$.
So there exists an irreducible polynomial $h\in\G_0[X]$ such that $\mathcal P=(h)$,
and $h$ is a nontrivial common divisor of the $f_i$ over $\G_0$, a contradiction.
 Thus,  $f_1,\ldots,f_\ell$ are  relatively prime over  $\G.$
\qed

\section{Partial Differential L\"{u}roth's Theorem}
For ordinary differential fields, Ritt and Kolchin proved the ordinary differential analog of  L\"uroth's theorem \cite{Ritt, Kolchin1944, Kolchin1947},
which states that every intermediate $\delta$-field $\G$ between $\ff$ and a   purely transcendental $\delta$-extension field $\ff\langle u\rangle$  
is always a simple extension.
Kolchin also pointed out in \cite[p. 164]{Kolchin1973} that his proof could be adapted to prove the generalized differential L\"uroth's theorem where 
$u$ is replaced by a finite set of $\delta$-indeterminates and $\G$ is of $\delta$-transcendence degree 1 over $\ff$.
But as Kolchin showed in \cite{Kolchin1944}, the L\"uroth's theorem does not always hold for partial differential fields. In this section, based on the classical proof of differential L\"uroth's theorem  \cite{Kolchin1947},
we aim to give a necessary and sufficient condition 
so that a partial differential analog of (generalized) L\"{u}roth's Theorem  holds.

\begin{theorem} \label{thm-Dluroth}
Let $\ff$ be a $\Delta$-field of characteristic 0 and $\{u_1, \ldots, u_n\}\subset\E$ a set of $\Delta$-indeterminates over $\ff$. 
If  $\G$ is a $\Delta$-field with $\ff \subsetneq \G \subseteq \ff\langle u_1, \ldots, u_n\rangle$ 
and the differential dimension polynomial of the point $\textbf{u}=(u_1, \ldots, u_n)$ over $\G$ is  of the form
\begin{equation}\label{eq-kolpoly} \omega_{\textbf{u}/\G}(t) = n\binom{t+m}{m}-\binom{t+m-s}{m}
\end{equation}\noindent for some $s\in\mathbb N$, then there exists  $v\in\G$ such that $\ff=\G\langle v\rangle$. 
Moreover, if $v_1$ and $v_2$ are two such elements, then there exist $a,b,c,d\in\ff$ such that 
$v_2=\frac{av_1+b}{cv_1+d}$.
\end{theorem} 

\proof   Let $\textbf{z}=(z_1,\ldots,z_n)$ be an $n$-tuple of  $\Delta$-variables over $\ff\langle u_1, \ldots, u_n\rangle$.
Let  $\Sigma=\mathbb I_{\G}(\textbf{u})$ be the vanishing ideal of  the point $\textbf{u}$  in $\G\{\textbf{z}\}$.
 Since the differential dimension polynomial of $\Sigma$ is of the form (\ref{eq-kolpoly}),
by Lemma \ref{lm-codim1},   there exists an irreducible $\Delta$-polynomial
 $A\in\G\{\textbf{z}\}$ such that $\Sigma=\sat(A)$.
The coefficients of $A$ are in $\G\subseteq \ff\langle \bu\rangle$, 
and thus are  quotients of elements in $\ff\{\bu\}$
and we may assume that one of its coefficients is $1$.

Let $v=\frac{H(\textbf{u})}{K(\textbf{u})}\in\G\backslash\ff$ be a coefficient of $A(\z)$ with $H, K\in\ff\{\bu\}$ coprime.
Then $$v K(\textbf{z})- H(\textbf{z})\in\Sigma.$$
In the following, we proceed to show that $\G=\ff\langle v\rangle$.
Let $D(\textbf{u})\in\ff\{\bu\}$ be the least common multiple of denominators of coefficients of $A$
and set $B(\textbf{u},\textbf{z})= D(\textbf{u})A(\textbf{z})$. 
Clearly, $B(\textbf{u},\textbf{u})=0.$
Fix $\rr$ to be an elimination ranking $\bu\prec\z$ on $\ff\{\bu, \z\}$ with  $\theta(y_j)<\tau(z_k)$ for all $\theta,\tau\in\Theta$. 
Without loss of generality, we may assume that $\lead(B)= \theta_0(z_n)$.

\vskip3pt
Let $B=\prod_{i=1}^{p}B_i$ be the irreducible factorization of $B$ in $\ff\{\bu,\z\}$. 
We first claim that \begin{enumerate}
\item[C1)] For each $i=1,\ldots,p$, $\lead(B_i)=\lead(B)=\theta_0(z_n)$;
\item[C2)] For $1\leq i\neq j\leq p$, $B_i\neq B_j.$
\end{enumerate}

To prove C1), first note that each $B_i$ effectively involves $\z$.
Indeed, if there is some $B_{i_0}\in\ff\{\bu\}$, then all the coefficients of $B$ regarded as a $\Delta$-polynomial in $\z$ are divisible by $B_{i_0}$.
In particular, since one of the coefficients of $A(\z)$ is $1$,  $D(\bu)$ is divisible by $B_{i_0}$.
Thus, $D(\bu)/B_{i_0}(\bu)$ is  a common multiple of the denominators of all the coefficients of $A(\z)$ with degree less than $D(\bu)$, 
a contradiction.
Now we show that each $B_i$ has the same leader as $B$, i.e., $\lead(B_i)=\theta_0(z_n)$. 
Rewrite $B(\bu, \z)$ as a univariate polynomial in  $\theta_0(z_n)$, then 
\begin{equation} B(y,z)=\sum_{j}\lambda_j(\bu,\z)(\theta_{0}z_n)^j
\end{equation} \noindent
with each $\lambda_j\in\ff\{\bu, \z\}$ free from $\theta_{0}z_n$.
 If for some $j_0$, $\lead(B_{j_0})\neq \theta_{0}z_n$, then
each $\lambda_j$ is divisible by   $B_{j_0}$.  
This implies that all the $\frac{\lambda_j}{D(\bu)}\in\G\{\z\}$ have a common divisor $\frac{B_{j_0}}{D(\bu)}$ in $\ff\langle\bu\rangle\{\z\}$.
By Lemma \ref{lm-coprime},  they also have a nontrivial common divisor in $\G\{\z\}$. 
Thus,   $A$ has a nontrivial divisor in $\G\{\z\}$, contradicting  the irreducibility of $A$.
So for each $B_i$, $\lead(B_i)=\theta_0(z_n)$.

\vskip3pt

 To prove C2), since $A(\z)$ is irreducible in $\G\{\z\}$, $A(\z)$ and $\S_A$ are coprime 
 in $\G\{z\}$,  and also coprime in $\ff\langle \bu\rangle\{\z\}$ by Lemma \ref{lm-coprime}. 
 So, in $\ff\langle \bu\rangle\{\z\}$, $\gcd(B,\S_B)=\gcd(D(\bu)A(\z),D(\bu)\S_A)=1$.
As each $B_j$ effectively involves $\theta_0(z_n)$ by C1),  
$B$ cannot have  multiple divisors in $\ff\{\bu,\z\}$ and  thus all the $B_i$ are distinct.

\vskip5pt
Assume all the derivatives in $\Theta(\bu)$ which appear  effectively in $B(\bu, \z)$ are 
$\theta_{1}(u_{e_1})<\theta_{2}(u_{e_2})< \cdots< \theta_{r}(u_{e_r})$. 
In the irreducible factorization of $B(\bu, \z)$,   we suppose without loss of generality that $B_1(\bu, \z)$ effectively involves $\theta_{r}(u_{e_r})$.
Then for each $j\neq1$, $B_j\not\in\sat(B_1)$,  for by the above claim, $\lead(B_j)=\lead(B_1)$ and $B_j$ can not be divisible by $B_1$.
Let $(\bu,\zeta_1)$ be a generic point of $\sat(B_1)\in\ff\{\bu,\z\}$. 
Note that $B(\bu,\z)\in \sat(B_1)$ and $\S_B\not\in \sat(B_1)$, for $\S_{B_1}\notin\sat(B_1)$. 
So $\zeta_1$ is a non-singular solution of $A(\z)\in\G\{\z\}$ and thus a solution of $\Sigma.$ 
Thus, $vK(\z)-H(\z)$ vanishes at $\zeta_1$. 
Therefore, 
$H(\bu)K(\z)-K(\bu)H(\z)\in\sat(B_1)$.
Consider temporarily the elimination ranking $\mathscr{R}_1: \z\prec\bu$ where the rankings on $\Theta(\z)$ and $\Theta(\bu)$ are the ones induced by $\mathscr{R}$.
Under $\mathscr{R}_1$,   $H(\bu)K(\z)-K(\bu)H(\z)$ is of leader  not higher than $\lead_{\mathscr{R}_1}(B_1) =\theta_{r}(u_{e_r})$ and thus   is partially reduced w.r.t. $B_1$.
By  Lemma \ref{le-generalcomponent}, $H(\bu)K(\z)-K(\bu)H(\z)$ is divisible by $B_1$.

Similarly, we can show that $H(\bu)K(\z)-K(\bu)H(\z)$ is divisible by all the $B_j$ which effectively involves $\theta_{r}u_{e_r}$, say, $B_1,\ldots, B_s$.
Thus, we get a relation
\begin{equation} H(\bu)K(\z)-K(\bu)H(\z)=L(\bu,\z)B_1(\bu,\z)\cdots B_s(\bu,\z)
\end{equation}
Note that the degree of $\prod_{i=1}^sB_i$ in $\theta_{r}(u_{e_r})$ is the same as $B(\bu, \z)$
and the degree of $H(\bu)K(\z)-K(\bu)H(\z)$ in any subset of $\Theta(\bu)$ is not greater than that of $B(\bu,\z)$. 
Thus $L(\bu, \z)$ does not involve $\theta_{r}(u_{e_r})$.

Let  $B_{s+1}$ be an irreducible factor effectively involving  $\theta_{r-1} (u_{e_{r-1}})$,
and $(\bu,\zeta_{s+1})$ be a generic point of $\sat(B_{s+1})$.
Similarly as above, $\zeta_{s+1}$ is a solution of $\Sigma$,
and thus, $H(\bu)K(\z)-K(\bu)H(\z)\in\sat(B_{s+1})$. 
Since $B_{j}\notin \sat(B_{s+1})$ for $j=1,\ldots,s$,  $L(\bu,\z)\in\sat(B_{s+1})$.
Under $\mathscr{R}_1$,   $L$ is of leader not higher than $\lead_{\mathscr{R}_1}(B_{s+1}) =\theta_{r-1}(u_{e_{r-1}})$ and thus  is partially reduced w.r.t. $B_{s+1}$.
By  Lemma \ref{le-generalcomponent}, $L$ is divisible by   $B_{s+1}$. 
Let $B_{s+1},\ldots ,B_{t}$ be all the irreducible factors of $B(\bu,\z)$ whose highest derivative in $\Theta(\bu)$ is $\theta_{r-1}u_{e_{r-1}}$.
Then $L(\bu,\z)=M(\bu,\z)B_{s+1}\cdots B_{t}$.
By comparing the degrees of $B(\bu,\z)$ and $M(\bu,\z)B_1\cdots B_{t}$ in $\theta_{r-1}(u_{e_{r-1}})$ and $\theta_{r}(u_{e_r})$, 
 $M(\bu,\z)$ does not involve $\theta_{r-1}(u_{e_{r-1}})$ and $\theta_{r}(u_{e_r})$. 

\vskip5pt
Continuing this process, we finally obtain a relation
\begin{equation}
 H(\bu)K(\z)-K(\bu)H(\z)=P(\z)B_1(\bu,\z)\ldots B_p(\bu,\z)=P(\z)B(\bu,\z)
\end{equation} \noindent
where $P(\z)\in\ff\{\z\}$.
Since $H(\z)$ and $K(\z)$ are coprime,  $H(\bu)K(\z)-K(\bu)H(\z)$ has no nontrivial divisor in $\ff\{\z\}$, and  $P(\z)=a\in\ff$ follows. 
Thus, $v K(\z)-H(\z)=\frac{aD(\bu)}{H(\bu)} \cdot A(\z)$.
Note that $\frac{aD(\bu)}{H(\bu)}$ is free from $\Theta(\z)$.
Thus, we have
\begin{equation}
v K(\z)-H(\z)=\alpha A(\z)  \text{ for some }  \alpha \in \G.
\end{equation}

\vskip3pt
We are ready to  show $\G=\ff\langle v\rangle$. 
Let $\gamma=\frac{U(\bu)}{V(\bu)}\in\G\backslash\ff$ with $\gcd(U,V)=1$.
Then $\gamma V(\z)-U(\z) \in \Sigma$. 
Let $\phi$ be any $\Delta$-isomorphism of $\ff\langle \bu\rangle$ w.r.t. $\ff\langle v\rangle$ and $\phi(v_i)=\bar{v}_i$ for $1\leq i\leq n$.
Then $$0=\phi(0)=\phi(v K(\bu)-H(\bu))=v  K(\bar{\bu})-H(\bar{\bu})=\alpha A(\bar{\bu}).$$
So $A(\bar{\bu})=0$.
Since $\S_A=\frac{\partial A}{\partial \theta_0(z_k)}$ does not vanish at $\bu$, $\frac{\partial (v K(\z)-H(\z))}{\partial \theta_0(z_k)}\big|_{\z=\bu}\neq0$.
So $\phi\big(\frac{\partial (v K(\z)-H(\z))}{\partial \theta_0(z_k)}\big|_{\z=\bu}\big)=\frac{\partial (v K(\z)-H(\z))}{\partial \theta_0(z_k)}\big|_{\z=\bar{\bu}}=\alpha\S_A(\bar{\bu})\neq0$.
Thus, $\bar{\bu}$ is a zero of $\Sigma=\sat(A)$.
Consequently, $\gamma V(\bar{\bu})-U(\bar{\bu})=0$, and
$\phi(\gamma)=\frac{U(\bar{\bu})}{V(\bar{\bu})}=\gamma$ follows. 
By Lemma \ref{Kolchin-isomorphism}, $\gamma\in\ff\langle v\rangle$ and finally  $\G=\ff\langle v\rangle$ is proved. 

\vskip3pt
 {It remains to show the uniqueness of L\"{u}roth generators under fractional linear transformations.
Suppose $v_1$ and $v_2$ are two L\"{u}roth generators of $\G/\ff$. 
Since $\ff\langle v_1\rangle=\ff\langle v_2\rangle$, there exist $U_1, V_1, U_2, V_2\in\ff\{z\}$ such that 
$v_2=\frac{U_1(v_1)}{V_1(v_1)}$ and $v_1=\frac{U_2(v_2)}{V_2(v_2)}$.
Consider the vanishing differential ideal $\mathcal I=\mathbb I_{\ff}(v_1, v_2)\subset\ff\{z_1,z_2\}$ of $(v_1,v_2)$ over $\ff.$
Le $\mathcal A=A_{11},\ldots,A_{1\ell_1}$ be a characteristic set of $ \mathcal I$ under the elimination ranking $z_1<z_2$.
Clearly, the leading variable of each $A_{1i}$ is $z_2$.
Since $H:=U_1(z_1)\cdot z_2 -V_1(z_1)\in  \mathcal I$ can be reduced to zero by $\mathcal A$ and  $\lead(H)=z_2$, 
$\ell_1=1$ and $\lead(A_{11})=z_2$.
Thus, $A_{11}=P(z_1)z_2+Q(z_1)$ and $\mathcal I=\sat(A_{11})$ follows.
By Lemma \ref{le-generalcomponent}, $A_{11}$ is also a characteristic set of $\mathcal I$ under the elimination ranking $\mathscr{R}_2: z_2<z_1$.
The fact that $V_2(z_2)z_1-U_2(z_2)\in\mathcal I$ implies that the leader of $A_{11}$ under $\mathscr{R}_2$ is $z_1$.
Thus, there exist $a, b, c, d\in\ff$ with $ad-bc\neq0$ s.t. $A_{11}=(az_2+b)z_1+(cz_2+d)$.
Hence, $v_1=-(cv_2+d)/(av_2+b)$.
}
\qed

\begin{remark}
 The element $v$ in Theorem \ref{thm-Dluroth} is called a {\rm L\"{u}roth generator} of $\G$.
If  $\G$ is a simple $\Delta$-extension, we say a {L\"{u}roth generator} exists for $\G$.
Theorem \ref{thm-Dluroth} gives a sufficient condition such that L\"{u}roth generators exist for an intermediate $\Delta$-field of $\ff\subseteq\ff\langle u_1,\ldots,u_n\rangle$,
and it claims  the uniqueness of L\"{u}roth generators  up to  linear fractional transformations.
 \end{remark}

\begin{example} \label{ex-luroth1}
Let $m=2$ and we consider the $\Delta$-field extensions $\ff \subsetneq \G \subseteq \ff\langle u\rangle$ with $\G=\ff\langle \delta_1\delta_2(u),\delta_1(u)+\delta_1\delta_2^2(u)\rangle$.
Note that $\omega_\G(t)={t+2\choose 2}-{t+1\choose 2}$,
so $\G$ admits a L\"{u}roth generator.
It is easy to verify that $\G=\ff\langle \delta_1(u)\rangle$.
\end{example}

Given an intermediate $\Delta$-field $\G$ between  $\ff $ and $\ff\langle u_1, \ldots, u_n\rangle$,
Theorem \ref{thm-Dluroth} shows that (\ref{eq-kolpoly}) is a sufficient condition for the existence of a L\"{u}roth generator.
We now show that (\ref{eq-kolpoly}) is also a necessary condition such that  a L\"{u}roth generator exists for $\G$.

\begin{cor} \label{cor-lurothcondition}
Let  $\G$ be a $\Delta$-field with $\ff \subsetneq \G \subseteq \ff\big\langle u_1, \ldots, u_n\big\rangle$ with $\ff$ and the $y_i$ as in Theorem \ref{thm-Dluroth}. Then a necessary and sufficient condition for $\G$ to have a L\"{u}roth generator is that
 the differential dimension polynomial of $\textbf{u}=(u_1, \ldots, u_n)$ over $\G$ is of the form
\begin{equation} \label{eq-lurothnece}  \omega_{\textbf{u}/\G}(t) = n\binom{t+m}{m}-\binom{t+m-s}{m}
\end{equation} for some $s\in\mathbb N$. 
\end{cor}
\proof     
By Theorem \ref{thm-Dluroth}, we only need to show the necessity part.
 Suppose $\G$ has a L\"ur{o}th generator $\frac{P(\textbf{u})}{Q(\textbf{u})}$ in reduced form, i.e., 
 $\G=\ff\langle\frac{P(\textbf{u})}{Q(\textbf{u})}\rangle$.
Take an orderly ranking  $\mathscr{R}$ on $\ff\{u_1,\ldots,u_n\}$ and suppose without loss of generality that $\max\{\lead(P),\lead(Q)\}=\theta(u_n)$. 
Let $s:=\ord(\theta)$ and  for  $t\geq s$, set $\G_{t}=\G(u_1^{[t]},\ldots,u_{n-1}^{[t]})$.
Then
\begin{eqnarray} \label{eq-kolcpolyofy}
\omega_{\textbf{u}/\G}(t)&=&\trdeg\,\G(u_1^{[t]},\ldots,u_n^{[t]})\Big/\G  \nonumber \\ 
& =&  \trdeg\,\G_{t}/\G+\trdeg\,\G_{t}(u_{n}^{[t]})/\G_{t}.   
\end{eqnarray} 
 Note that $u_1^{[t]},\ldots,u_{n-1}^{[t]}$ are algebraically independent over $\G=\ff\langle\frac{P(\textbf{u})}{Q(\textbf{u})}\rangle.$ For otherwise, there would be some $\ell$ s.t. $u_1^{[t]},\ldots,u_{n-1}^{[t]}$ are algebraically dependent over $\ff\Big(\big(\frac{P(\textbf{u})}{Q(\textbf{u})}\big)^{[\ell]}\Big).$
Consequently, $\big(\frac{P(\textbf{u})}{Q(\textbf{u})}\big)^{[\ell]}$ are algebraically dependent over $\ff(u_1^{[t]},\ldots,u_{n-1}^{[t]})$, contradicting  the fact that the Jacobian matrix of $\big(\frac{P(\textbf{u})}{Q(\textbf{u})}\big)^{[\ell]}$ w.r.t. $\Theta_{\leq \ell}(\theta(y_n))$ is of full rank \cite[p. 12, 11.4]{Lefschetz}. Thus, 
$\trdeg\,\G_{t}/\G=(n-1)\binom{t+m}{m}.$

On the other hand, setting $\ff'=\ff\langle u_1,\ldots,u_{n-1}\rangle$,
 by Lemma \ref{lm-ukolchinpol}, we have $\omega_{u_n/\ff'\langle\frac{P(\textbf{u})}{Q(\textbf{u})}}(t)=\trdeg\,\ff'\langle\frac{P(\textbf{u})}{Q(\textbf{u})}\rangle(u_{n}^{[t]})/\ff'\langle\frac{P(\textbf{u})}{Q(\textbf{u})}\rangle=\binom{t+m}{m}-\binom{t+m-s}{m}.$
 As $\G_t\subset\ff'\langle\frac{P(\textbf{u})}{Q(\textbf{u})}\rangle$,
 $\trdeg\,\G_{t}(u_{n}^{[t]})/\G_{t}\geq \binom{t+m}{m}-\binom{t+m-s}{m}$.
Similarly as in the proof of Lemma \ref{lm-ukolchinpol}, we can show  that 
$$\Theta_{\leq t-s}\theta(u_n)\subset\ff\big(u_1^{[t]},\ldots,u_{n-1}^{[t]} ,u_n^{[t]}\backslash\Theta_{\leq t-s}\theta(u_n), \theta(u_n)\big).$$
So $\trdeg\,\G_{t}(u_{n}^{[t]})/\G_{t}\leq\trdeg\,\G_t\big(u_n^{[t]}\backslash\Theta_{\leq t-s}\theta(u_n), \theta(u_n)\big)/\G_t
=\binom{t+m}{m}-\binom{t+m-s}{m}.$
Thus,  by (\ref{eq-kolcpolyofy}), $\omega_{\textbf{u}/\G}(t) = n\binom{t+m}{m}-\binom{t+m-s}{m}$ follows.
\qed

\vskip5pt
The following example was given by Kolchin in \cite[Example 1]{Kolchin1944} where he proved the given partial differential intermidiate field does not admit a 
L\"{u}roth generator.
Corollary \ref{cor-lurothcondition} provides an alternative simple proof. 
\begin{example}
Let $m=2$ and we consider the $\Delta$-field extensions $\ff \subsetneq \G \subseteq \ff\langle u\rangle$ with $\G=\ff\langle \delta_1(u),\delta_2(u)\rangle$.
Note that $\omega_{u/\G}(t)=1\neq {t+2\choose 2}-{t+2-s\choose 2}$ for any $s\in\mathbb N$,
so by Corollary \ref{cor-lurothcondition}, $\G$ does not admit a L\"{u}roth generator. 
\end{example}

\begin{remark}
In the ordinary differential case (i.e., $m=1$),  Corollary \ref{cor-lurothcondition} implies that 
$\G$ is a simple differential extension of $\ff$ if and only if $\omega_{\textbf{u}/\G}(t) = (n-1)(t+1)+s$ for some $s\in\mathbb N$. 
The latter condition is equivalent to that $\G$ is of differential transcendence degree one over $\ff$.
This is indeed the generalized differential L\"uroth's theorem in the case $n>1$ (see \cite[p. 164]{Kolchin1973}  and \cite{DJS2018}) and the classical differential L\"uroth's theorem \cite{Ritt1932, Kolchin1944} in the case $n=1$.
\end{remark}

\section{On the algorithmic aspects  of partial differential L\"{u}roth's theorem}

\indent As we have shown in the last section, Corollary \ref{cor-lurothcondition} gives a criterion to determine whether a  given   intermediate  $\Delta$-field $\G$ between  $\ff$ and $\ff\langle u_1, \ldots, u_n\rangle$ admits a L\"{u}roth generator, and in the affirmative case, the proof of Theorem \ref{thm-Dluroth} gives a theoretical method to compute a L\"{u}roth generator. 
In this section, we focus on the algorithmic aspects of the partial differential L\"{u}roth theorem.

\medskip
Let $\G=\ff\langle \frac{P_1(\bu)}{Q_1(\bu)}, \,\ldots, \,\frac{P_r(\bu)}{Q_r(\bu)}\rangle$ be a $\Delta$-subfield of $\ff\langle \bu\rangle$,
where $\bu=(u_1, \ldots,u_n)$, and for each $j$, $\frac{P_j(\bu)}{Q_j(\bu)}\notin\ff$  is assumed to be  in reduced form. 
Let $$PS:=\{ x_1Q_1(\bu)-P_1(\bu),\ldots, x_rQ_r(\bu)-P_r(\bu) \}$$
and $$DS:=\{Q_1(\bu),\,\ldots,\,Q_r(\bu)\}.$$
Set $$\mathcal P=[PS]:(\prod_{j=1}^rQ_j)^\infty\subseteq\ff\{u_1,\ldots,u_n,x_1,\ldots,x_r\}$$ 
and a point
\begin{equation} \label{eq-eta}
\eta=\Big(u_1,\ldots,u_n, \frac{P_1(\bu)}{Q_1(\bu)}, \ldots, \frac{P_r(\bu)}{Q_r(\bu)}\Big)\in\ff\langle \bu\rangle^{n+r}.
\end{equation}
It is easy to show that $\mathcal P$ is a prime $\Delta$-ideal with a generic point $\eta$ and  $PS$ is a characteristic set of $\mathcal P$ under any elimination ranking\footnote{Unlike the ordinary differential case, there are  distinct choices of elimination rankings with the same variable ordering $x_1<x_2$. 
For example, apart from  $\theta(x_1)>\phi(x_2)$ for any $\theta, \phi$, we can choose $\theta(x_1)>\phi(x_1)$ if $\theta>\phi$   for the lexicographic (graded lexicographic, graded reverse  lexicographic) order given by $\delta_1>\cdots>\delta_m$ (same or different choice for $\Theta(x_2)$).
For convenience, in this paper, by saying {\it the elimination ranking $x_1\prec\cdots\prec x_r$}, our default ranking on each variable  is the canonical ranking, 
i.e., $\prod_{i=1}^m\delta_i^{k_i}(x_j) <\prod_{i=1}^m\delta_i^{l_i}(x_j)$ if and only if $ (\sum_{i} k_i, k_1,\ldots,k_{m})<_{lex}(\sum_{i} l_i, l_1,\ldots,l_{m})$.  } $\mathscr{R}_1: u_1\prec\cdots\prec u_n\prec x_1\prec\cdots\prec x_r$.

\medskip
Take another elimination ranking $\mathscr{R}_2$ with $ x_1\prec\cdots\prec x_r\prec u_1\prec\cdots\prec u_n$.
We perform the Wu-Ritt's zero decomposition algorithm as in Theorem \ref{th-zerodecom} for $\mathbb V(PS/DS)$ under the ranking $\mathscr{R}_2$, and we will get a finite number of $\Delta$-irreducible and coherent autoreduced sets $\mathcal C_1,\ldots,\mathcal C_l$ such that $$\mathbb V(PS/DS)=\bigcup_{i=1}^l\mathbb V(\sat(\mathcal C_i)/DS).$$

\begin{lemma} \label{lm-ck}
There exists a unique $\mathcal C_k$ such that $\eta$ is a zero of $\mathcal C_k$.
And for this $\mathcal C_k$, we have $\mathcal{P}=\sat(\mathcal C_k)$.
\end{lemma}
\proof  Note that $\eta$ is a point of $\mathbb V(PS/DS)$.
So there exists some $k$ such that $\eta\in\mathbb V(\sat(\mathcal C_k)/DS)\neq\emptyset$.
Since $\eta$ is a generic point of $\mathcal P$, $\sat(\mathcal C_k)\subseteq \mathcal P$.
Conversely, if $\xi$ is a generic point of $\sat(\mathcal C_k)$,
$\xi\in\mathbb V(\sat(\mathcal C_k)/DS)\subset \mathbb V(PS/DS)$.
So $\xi$ is a zero of $\mathcal P$ and consequently, $\mathcal P\subseteq\sat(\mathcal C_k)$ follows.
Thus, for such a $\mathcal C_k$,  $\mathcal{P}=\sat(\mathcal C_k)$,
and the uniqueness of $\mathcal C_k$ follows.
\qed 

\medskip
By Lemma \ref{lm-ck}, $\mathcal C_k$ is a characteristic set of $\mathcal P$ under the ranking $\mathscr{R}_2$.
Note that if $U:=\{x_{i_1},\ldots,x_{i_d}, u_{\sigma_1},\ldots,u_{\sigma_e}\}$ is a parametric set of $\mathcal C_k$,
then $\mathcal C_k$ is also a $\Delta$-irreducible and coherent autoreduced set under the elimination ranking $x_{i_1}\prec\cdots\prec x_{i_d}\prec\cdots\prec x_{i_r}\prec u_{\sigma_1}\prec\cdots u_{\sigma_e}\prec \cdots\prec u_{\sigma_n}$, where the rankings on $\Theta(x_i)$ and $\Theta(u_j)$
are the ones induced by the original ranking $\mathscr{R}_2$.
Also note that the differential dimension of $\mathcal P=\sat(\mathcal C_k)$ is $n$.
After changing the names of the variables, we may write $\mathcal C_k$ as follows:
\begin{equation}\label{eq-charset} \begin{array}{l} 
A_{1,1}(x_1,\ldots,x_{d},x_{d+1}), \ldots, A_{1,o_1}(x_1,\ldots,x_{d},x_{d+1}) \\
\cdots\cdots\\
A_{r-d,1}(x_1,\ldots,x_{d}, \ldots,x_r), \ldots, A_{r-d,o_{r-d}}(x_1,\ldots,x_{d}, \ldots,x_r) \\
B_{1,1}(x_1,\ldots,x_r,u_1,\ldots, u_{n-d+1}), \ldots, B_{1,l_1}(x_1,\ldots,x_r,u_1,\ldots, u_{n-d+1})\\
\cdots\cdots\\
B_{d,1}(x_1,\ldots,x_r,u_1,\ldots,u_{n}), \ldots, B_{d,l_d}(x_1,\ldots,x_r,u_1,\ldots,u_{n})
\end{array}
\end{equation}
where $x_{d+i}$ is the leading variable of each $A_{i,j}$, 
and $u_{n-d+i}$ is the leading variable of each $B_{i,j}$.
We now give a criterion for the existence of a L\"{u}roth generator for $\G$ based on the obtained $\mathcal C_k$ in the form   (\ref{eq-charset}).

\begin{theorem} \label{th-lurothcriterion}
Let  $\G=\ff\big\langle\frac{P_1(\bu)}{Q_1(\bu)}, \ldots, \frac{P_r(\bu)}{Q_r(\bu)}\big\rangle\subseteq\ff\langle \bu\rangle$.
A necessary and sufficient condition for $\G$ to have a L\"{u}roth generator is that  $d=1$ and $l_1=1$.
\end{theorem}
\proof  Introduce new variables $z_1,\ldots, z_n$ which are $\Delta$-independent over $\ff\langle\bu\rangle.$ For each $B_{i,j}\,(i=1,\ldots,d; j=1,\ldots,l_i)$ in (\ref{eq-charset}), let 
$$\bar{B}_{i,j}=B_{i,j}\Big(\frac{P_1(\bu)}{Q_1(\bu)},\ldots,\frac{P_r(\bu)}{Q_r(\bu)},z_1,\ldots, z_{n-d+i}\Big).$$
Then $\bar{B}_{i,j}\in\G\{z_1,\ldots,z_n\}$.
Take the elimination ranking $z_1\prec\cdots\prec z_n$ where the ranking on $\Theta(z_i)$ is  the one induced by  $\mathscr{R}_2$  on $\Theta(u_i)$.
Then $\bar{B}_{i,j}$ is a $\Delta$-irreducible and coherent autoreduced set in $\G\{z_1,\ldots,z_n\}$.
Denote the vanishing $\Delta$-ideal of $\bu$ in $\G\{z_1,\ldots,z_n\}$ by $\mathbb I_{\G}(\bu)$.
We now show that 
\begin{equation} \label{eq-inducedcharset}
\mathbb I_{\G}(\bu)=\sat\big((\bar{B}_{i,j})_{ i=1,\ldots,d; j=1,\ldots,l_i}\big). 
\end{equation}

Since $\eta$ given in (\ref{eq-eta}) is a generic point of $\mathcal P$ and (\ref{eq-charset}) is a characteristic set of $\mathcal P$, $\bar{B}_{i,j}(\bu)=0$ and $H_{\bar{B}_{i,j}}(\bu)=H_{B_{i,j}}(\eta)\neq0$.
Thus, $\sat\big((\bar{B}_{i,j})_{ i=1,\ldots,d; j=1,\ldots,l_i}\big)\subseteq\mathbb I_{\G}(\bu).$
On the other hand, suppose $f$ is an arbitrary polynomial in $\mathbb I_{\G}(\bu).$  
By clearing  denominators when necessary, we may suppose 
 there is an $F\in\ff\{x_1,\ldots,x_r,u_1,\ldots,u_n\}$
and $f=F(\frac{P_1(\bu)}{Q_1(\bu)},\ldots,\frac{P_r(\bu)}{Q_r(\bu)},z_1,\ldots,z_n)$. 
Since $f\in\mathbb I_{\G}(\bu)$, $F\in\sat(\mathcal C_k)$.
So there exist  $a\in\mathbb N$ such that $$\Big(\prod_{i,j}\prod_{p,q}\I_{A_{p,q}}\S_{A_{p,q}}\I_{B_{i,j}}\S_{B_{i,j}}\Big)^a\cdot F\in [A_{p,q}, B_{i,j}].$$
Substituting $x_i=\frac{P_i(\bu)}{Q_i(\bu)}$  and $u_j=z_j$ for all $i,j$ in the above identity, we have
 $$\Big(\prod_{i,j}\I_{\bar{B}_{i,j}}\S_{\bar{B}_{i,j}}\Big)^a\cdot f\in [\bar{B}_{i,j}: i=1,\ldots,d; j=1,\ldots,l_i].$$
Thus, $f\in\sat\big((\bar{B}_{i,j})_{ i=1,\ldots,d; j=1,\ldots,l_i}\big)$ and (\ref{eq-inducedcharset}) is proved.

\medskip
 By Corollary \ref{cor-lurothcondition}, a L\"{u}roth generator exists for $\G$ if and only if $\mathbb I_{\G}(\bu)$ is the general component of some $\Delta$-polynomial,
and by (\ref{eq-inducedcharset}), the latter holds if and only if $d=1$ and $l_1=1$.
\qed

\medskip
Based on Theorem \ref{th-lurothcriterion} and the proof of Theorem \ref{thm-Dluroth}, we now give an algorithm to decide 
whether a given finitely generated $\Delta$-extension field of $\ff$ contained in $\ff\langle \bu\rangle$ has a L\"{u}roth generator, and in the affirmative case, compute a L\"{u}roth generator.

\bigskip
\noindent\textbf{Algorithm Partial-Diff-L\"{u}roth}
\vskip3pt
\noindent\textbf{Input}: A finitely generated $\Delta$-field  $\G:=\ff \langle \frac{P_1(\bu)}{Q_1(\bu)}, $ $\ldots, \frac{P_r(\bu)}{Q_r(\bu)} \rangle\subseteq\ff\langle \bu\rangle$. 
\vskip3pt
\noindent\textbf{Output}: A L\"{u}roth generator $v$ of $\G$, if a L\"{u}roth generator exists for $\G$; \vskip2pt\hskip1truecm Otherwise, return ``$\G$ does not admit a L\"{u}roth generator".

\begin{itemize}
\item[1.]  Apply Wu-Ritt's zero decomposition theorem to $PS:=\{ x_1Q_1(\bu)-P_1(\bu),$ $\ldots, x_rQ_r(\bu)-P_r(\bu)\}$ and $DS:=\{Q_1(\bu),\ldots,Q_r(\bu)\} $ 
under the elimination ranking $x_1\prec\cdots\prec x_r\prec u_1\prec\cdots\prec u_n$  to compute an irreducible decomposition of $\mathbb V(PS/DS)$: 
$$\mathbb V(PS/DS)=\bigcup_{i=1}^p\mathbb V\big(\sat(\mathcal C_i)/DS),$$
where the obtained $\mathcal C_i$ are $\Delta$-irreducible and coherent autoreduced sets.

\item[2.] For $i=1,\ldots,p$, if $\mathcal C_i$ vanishes at $\eta=\big(\frac{P_1(\bu)}{Q_1(\bu)},  \ldots, \frac{P_r(\bu)}{Q_r(\bu)},u_1,\ldots,u_n\big)$, 
then set $\mathcal A:=\mathcal C_k$ and goto step 3; else, $i:=i+1$.

\item[3.] Let $T:=\{A\in\mathcal A\,|\,\text{lv}(A)=u_j \text{ for some } j\}$.
If $\text{Card}(T)\neq 1$, then return ``$\G$ does not admit a L\"{u}roth generator";
else, write $A\in T$ in the form $A=\sum_{i=1}^qa_{M_i}(x_1,\ldots,x_r)M_i(\bu)$ with the $M_i$   distinct $\Delta$-monomials in $\bu.$ 

\item[4.]  For $i=2,\ldots,q$, compute $v_i=\frac{a_{M_i}}{a_{M_1}}\Big|_{x_j=\frac{P_j(\bu)}{Q_j(\bu)},\, j=1,\ldots,r}$, and if $v_i\notin\ff$, then return $v_i$; else $i=i+1$.
 
\end{itemize}

\textbf{Proof of Correctness:}
\vskip2pt
For step 2, Lemma \ref{lm-ck} guarantees there is a unique $\mathcal C_k$ vanishing at $\eta$. 
For the obtained $\A$, by Theorem \ref{th-lurothcriterion}, 
a L\"{u}roth generator exists for $\G$ if and only if there is only one differential polynomial in $T$ with leading variable belonging to $\{u_1,\ldots,u_n\}$, i.e., $\text{Card}(T)=1$.
If $ T=\{A\}$, by (\ref{eq-inducedcharset}), $\bar{A}=A(\frac{P_1(\bu)}{Q_1(\bu)},  \ldots, \frac{P_r(\bu)}{Q_r(\bu)}, z_1,\ldots,z_n)\in\mathcal{G}\{z_1,\ldots,z_n\}$
can serve as a characteristic set of $\mathbb I_{\G}(\bu)$.
Let $b=a_{M_1}\Big|_{x_j=\frac{P_j(\bu)}{Q_j(\bu)},j=1,\ldots,r}$ and set $B=\frac{1}{b}\bar{A}$. Note that $B$ has a coefficient 1 and is also a characteristic set of $\mathbb I_{\G}(\bu)$. Since $u_1,\ldots,u_n$ are $\Delta$-independent over $\ff$, 
$B$ has at least a coefficient not contained in $\ff$.
By the proof of Theorem \ref{thm-Dluroth}, such a coefficient $v\notin\ff$ of $\bar{A}$ is a L\"{u}roth generator of $\G$.

\begin{example} \label{ex-lurothalg}
Continue from Example \ref{ex-luroth1}.
In this example,  $m=2$ and $\G=\ff\langle \delta_1\delta_2(u),\delta_1(u)+\delta_1\delta_2^2(u)\rangle$.
Let $PS:=\{x_1-\delta_1\delta_2(u), x_2-\delta_1(u)-\delta_1\delta_2^2(u)\}$.
Using Wu-Ritt's zero decomposition theorem to $PS$ w.r.t. the elimination ranking $\mathscr{R}: x_1\prec x_2\prec u$,
 we get $\V(PS)=\V(\A)$ with $\A=-\delta_2(x_2)+\delta_2^2(x_1)+x_1, \delta_1(u)+\delta_2(x_1)-x_2.$
Since $T=\{ \delta_1(u)+\delta_2(x_1)-x_2\}$ consists of only with polynomial,
$\G$ admits a L\"{u}roth generator and $v=\delta_2(\delta_1\delta_2(u))-(\delta_1(u)+\delta_1\delta_2^2(u)) =-\delta_1(u)$ is a L\"{u}roth generator. Thus, $\G=\ff\langle \delta_1(u)\rangle$.
\end{example}

\begin{example}
Let $m=2$ and consider the $\Delta$-field extensions $\ff \subsetneq \G \subseteq \ff\langle u\rangle$ with $\G=\ff\langle \frac{\delta_1(u)}{u}, u+\delta_1(u), \delta_2(u)\rangle$.
Let $PS:=\{ux_1-\delta_1(u), x_2-u-\delta_1(u), x_3-\delta_2(u)\}$.
Using Wu-Ritt's zero decomposition theorem to $PS$ w.r.t. the elimination ranking $\mathscr{R}: x_1\prec x_2\prec x_3\prec u$,
 we get $\V(PS/u)=\V(\C_1/u\H_{\C_1})\cup\V(\C_2/u)$ with 
$\C_1=(x_1+1)\delta_1(x_2) - \delta_1(x_1)x_2 - x_1^2x_2 - x_1x_2, (x_1^2+2x_1+1)x_3-\delta_2(x_2)x_1- \delta_2(x_2) + \delta_2(x_1)x_2, (x_1+1)u-x_2$ and $\C_2=x_1 +1, x_2, \delta_1(x_3) + x_3, \delta_1(u) + u, \delta_2(u)-x_3.$
Only $\A_1$ vanishes at $\eta=(\delta_1(u)/u, u+\delta_1(u),\delta_2(u))$.
Set $\A=\C_1$ and then $T=\{ (x_1+1)u-x_2\}$ has cardinality 1. 
Thus,  $v=\frac{x_2}{x_1+1}|_{x=\eta}=u$ is a L\"{u}roth generator and $\G=\ff\langle u\rangle$.
\end{example}

\begin{example}
Let $m=2$ and $\G=\ff\langle \delta_1(u),\delta_2(u)\rangle\subseteq \ff\langle u\rangle$.
Let $PS:=\{x_1-\delta_1(u), x_2-\delta_2(u)\}.$
Using Wu-Ritt's zero decomposition theorem to $PS$ w.r.t.  the elimination ranking $\mathscr{R}: x_1\prec x_2\prec u$,  we get $\V(PS)=\V(\A)$ with $\A=\delta_1(x_2)-\delta_2(x_1), -\delta_1(u)+x_1, -\delta_2(u)+x_2.$
Since $T=\{ -\delta_1(u)+x_1, -\delta_2(u)+x_2\}$ has cardinality 2, $\G$ does not admit a L\"{u}roth generator.
\end{example}

\begin{remark}
In  {Algorithm Partial-Diff-L\"{u}roth}, instead of applying  Wu-Ritt's zero decomposition theorem to $PS$ and $DS$, we can perform Rosenfeld-Gr\"obner algorithm \cite{BOLP09, Hubert} to the system $PS=0$ and $DS\neq0$, which is factorization-free, and thus more efficient and has been implemented in Maple.
\end{remark}

We conclude this section by giving an optimal order bound for   L\"uroth generators in terms of the orders of  the given generators of a finitely generated $\Delta$-field extensions. 
The ideas come from \cite{DJS2014}.
 
Let   $v=\frac{f(\bu)}{g(\bu)}\in\ff\langle\bu\rangle$ be in reduced form.
For some   ranking $\mathscr{R}$ on $\Theta(\bu)$, the leader of $v$ under $\mathscr{R}$ is defined as 
$\lead_{\mathscr{R}}(v):=\max_{\mathscr{R}}\{\lead(f), \lead(g)\}$.
And if $\mathscr{R}$ is an orderly ranking, $\lead_{\mathscr{R}}(v)$ is called a {\it supporting leader} of $v$.
We denote the set of all {supporting leaders} of $v$ by $\text{Sld}(v)$.

\begin{theorem}
Let  $\G=\ff\big\langle\frac{P_1(\bu)}{Q_1(\bu)}, \ldots, \frac{P_r(\bu)}{Q_r(\bu)}\big\rangle\subseteq\ff\langle \bu\rangle$ be a finitely generated $\Delta$-extension of $\ff$. 
Suppose  $v$ is a L\"uroth generator  of $\G$.
Then we have
\begin{equation} \label{eq-Sld}
\text{\rm Sld}\big(\frac{P_i(\bu)}{Q_i(\bu)}\big)\subseteq\Theta\big(\text{\rm Sld}(v)\big),\, i=1,\ldots,r.
\end{equation}
In particular, the order of $v$ cannot exceed $ \min_i\{\ord\big(\frac{P_i(\bu)}{Q_i(\bu)}\big)\}.$ 
\end{theorem}
\proof  For each $i$, denote $v_i=\frac{P_i(\bu)}{Q_i(\bu)}$.
Then there exists coprime $\Delta$-polynomials $f_{i}, g_{i}\in\ff\{z\}$ s.t. $v_i=\frac{f_{i}(v)}{g_{i}(v)}$.
Fix an orderly ranking $\mathscr{R}$ on $\Theta(\bu)$.
 Let $\theta(u_k)=\lead_{\mathscr{R}}(v)$ and
  $\theta(u_{o_i})=\lead_{\mathscr{R}}(v_i).$ 
It suffices to show that $\theta(u_{o_i})\in\Theta(\theta(u_k))$.

\vskip4pt
By the claim ($\star$) in the proof of Lemma \ref{lm-ukolchinpol}, for each $\pi\in\Theta_{>0}$,
 $\rank(\pi(v))=(\pi(\theta(u_k)),1)$, so $\frac{\partial \pi(v)}{\partial \pi(\theta(u_k))}\neq0$.
For $i=1,\ldots,r$, let $$\pi_i(\theta(u_k)):=\max\big\{\pi(\theta(u_k))| \pi(z) \text{ effectively appears in } {f_{i}(z)}/{g_{i}(z)}\big\}.$$
Note that if we define a ranking on $\Theta(z)$ by $\tau_1(z)<\tau_2(z)$ if and only if $\tau_1(u_k)<_{\mathscr{R}}\tau_2(u_k)$,
then $\lead(f_i(z)/g_i(z))=\pi_i(z)$.
 Clearly, $\theta(u_{o_i})\leq \pi_i(\theta(u_{k}))$.
Since $\frac{\partial v_i}{\partial \pi_i(\theta(u_{k}))}=\sum_{\pi\in\Theta}\frac{\partial f_i(z)/g_i(z)}{\partial \pi(z)}(v)\cdot\frac{\partial \pi(v)}{\partial \pi_i(\theta(u_{k}))}=\frac{\partial f_i(z)/g_i(z)}{\partial \pi_i(z)}(v)\cdot\frac{\partial \pi_i(v)}{\partial \pi_i(\theta(u_{k}))}\neq0$,
$\theta(u_{o_i})=\pi_i(\theta(u_k))$.
Thus, (\ref{eq-Sld}) holds.
 \qed

\begin{example} Continue from Example \ref{ex-lurothalg}.
In this example,  $\text{\rm Sld}(v)=\{\delta_1(u)\}$, $\text{\rm Sld}(\delta_1\delta_2(u))=\{\delta_1\delta_2(u)\}\subseteq\Theta\big(\text{\rm Sld}(v)\big),  \text{\rm Sld}(\delta_1(u)+\delta_1\delta_2^2(u))=\{\delta_1\delta_2^2(u)\}\subseteq\Theta\big(\text{\rm Sld}(v)\big)$. 
\end{example}

\section{Applications to proper re-parametrization}

\noindent An (irreducible) differential curve is an (irreducible) differential variety $V\subset \mathbb A^n$ 
of differential dimension 1.
An irreducible differential curve $V$ is unirational if  there exists a dominant  differential  rational map from $\mathbb A^1$ to $V$, that is, there exist $\mathcal P(u)\in\ff\langle u\rangle^n$ s.t. the differential rational map $\phi: \mathbb A^1 { \dashrightarrow} V$ with  $\phi(a)=\mathcal P(a)$ is dominant (i.e., $V$ is the Kolchin closure of the image of $\phi$.) If additionally $\phi$ has an inverse dominant  differential  rational map,
then $V$ is called rational. In other words,  a differential curve  is  rational if it is differentially birationally equivalent to  $\mathbb A^1$.

A natural question is to ask whether a unirational differential curve is rational.
In the ordinary differential case, the differential L\"{u}roth's theorem  gives an affirmative answer.
But in the partial differential case, a unitatonal differential curve may not  be rational.
In this section, based on the partial differential L\"{u}roth's theorem,
 we will give a rationality criterion for partial differential curves. 
This is closely related to the proper re-parametrization problem as treated in \cite[Theorem 6.2]{Gao2003} for ordinary differential curves.
 

\medskip
 Given a set of $\Delta$-polynomials $P_i, Q_i\in\ff\{u_1\}$ for $ i=1,\ldots,n$ 
 with $\gcd(P_i,Q_i)=1$
 and not all $P_i/Q_i\in\ff$, we call 
 \begin{equation}\label{eq-defpdrpes}
 x_1=\frac{P_1(u_1)}{Q_1(u_1)},\,\,\ldots,\,\, x_n=\frac{P_n(u_1)}{Q_n(u_1)}.
 \end{equation}
a set of {\it partial differential rational parametrization equations} (PDRPEs) in the parameter $u_1$.
For ease of notation, denote $\mathcal P(u_1)=\big(\frac{P_1(u_1)}{Q_1(u_1)},\ldots, \frac{P_n(u_1)}{Q_n(u_1)}\big)$
and $\mathbb I_{\ff}(\mathcal P(u_1))=\big\{g\in\ff\{x_1,\ldots,x_n\}\big| \,g(\mathcal P(u_1))=0\big\}$. 
The $\Delta$-variety $V=\mathbb V(\mathbb I_{\ff}(\mathcal P(u_1)))$ is called the {\it implicit variety} of (\ref{eq-defpdrpes}) and accordingly, (\ref{eq-defpdrpes}) or  $\mathcal P(u_1)$ for simplicity  is called a {\it $\Delta$-rational parametrization} of $V$.
Since at least a $P_i/Q_i\notin\ff$, the implicit variety $V$ is a $\Delta$-curve which is unirational ($\mathcal P(u_1)$ defines a dominant $\Delta$-rational map from $\mathbb A^1$ to $V$.)

Li introduced the definition of proper PDRPEs  in the geometrical language in \cite{Li2006} following \cite{Gao2003}. We give an equivalent definition in the algebraic setting.
 \begin{definition}
 A set of PDRPEs (\ref{eq-defpdrpes}) is said to be {\rm proper}  if $$\ff\big\langle\frac{P_1(u_1)}{Q_1(u_1)},\ldots, \frac{P_r(u_1)}{Q_r(u_1)}\big\rangle=\ff\langle u_1\rangle.$$
 \end{definition}

%

\begin{lemma}
Let $V$ be a unirational $\Delta$-curve.
Then $V$ is rational if and only if $V$ has a proper $\Delta$-rational parametrization.
\end{lemma}
\proof First, suppose $V$ is rational. Then there exist dominant $\Delta$-rational maps 
\[ \begin{array}{ccc} \phi: \mathbb A^1 & { \dashrightarrow} & V\\ 
\quad u_1 &  {\longmapsto } &  \mathcal P(u_1)= \big(\frac{P_1(u_1)}{Q_1(u_1)},\ldots, \frac{P_n(u_1)}{Q_n(u_1)}\big)  \end{array} \]
and 
\[ \begin{array}{ccc} \varphi: V &{ \dashrightarrow}& \mathbb A^1\\ 
\qquad (x_1,\ldots,x_n) & \longmapsto &  U(x_1,\ldots,x_n)  \end{array} \]
which are inverse to each other.
Since $\phi$ is dominant and $\varphi\circ\phi(u_1)=U(\mathcal P(u_1))=u_1$, $ \mathcal P(u_1)$ is a proper $\Delta$-rational parametrization of $V$.

Conversely, suppose $V$ has a proper $\Delta$-rational parametrization $ \mathcal P(u_1)$.
Then $\ff\langle \mathcal P(u_1)\rangle=\ff\langle u_1\rangle$.
So there exists $U(x_1,\ldots,x_n)\in\ff\langle x_1,\ldots,x_n\rangle$   such that $u_1=U\big( \mathcal P(u_1)\big).$ Clearly, $\mathcal P(u_1)$ and $U(x_1,\ldots,x_n)$ define $\Delta$-birationally equivalent map between $\mathbb A^1$ and $V$.
\qed

 \medskip

Given a set of PDRPEs  of the form (\ref{eq-defpdrpes}),
 \cite[Theorem 6.1]{Li2006} provides a method to detect whether (\ref{eq-defpdrpes}) is proper.
But if  (\ref{eq-defpdrpes}) is not proper, the problem to decide whether (\ref{eq-defpdrpes})  has a proper re-parameterization is unsolved.
In the following, we will solve this problem via using  {Algorithm Partial-Diff-L\"{u}roth}.

\vskip3pt
Let $\mathcal S=\{x_1Q_1(u_1)-P_1(u_1), \ldots, x_nQ_n(u_1)-P_n(u_1)\}$. 
Under the elimination ranking $x_1\prec x_2\prec\cdots\prec x_r\prec u_1$,
$\sat(\mathcal S)$ has a characteristic set 
\begin{equation} \label{eq-charset5.3}\mathcal A:=
\left\{\begin{array}{l}
A_{1,1}(x_1,x_2),\ldots, A_{1,r_1}(x_1,x_2) \\
 
\cdots\cdots\cdots \\  
A_{n-11,1}(x_1,\ldots,x_n),\ldots, A_{n-1,r_{n-1}}(x_1,\ldots,x_n) \\
\\
B_{1,1}(x_1,\ldots,x_n,u_1),\ldots, B_{1,l_1}(x_1,\ldots,x_n,u_1). 
\end{array}
\right.
\end{equation}
  \cite[Theorem 6.1]{Li2006} tells that  (\ref{eq-defpdrpes}) is proper if and only if $l_1=1$ and $ B_{1,1}=I_1(x_1,\ldots,x_n)u_1-U_1(x_1,\ldots,x_n)$.
Suppose  (\ref{eq-defpdrpes}) is not proper.
By saying (\ref{eq-defpdrpes}) has a proper re-parameterization, we mean there exist  a new parameter $v=R(u_1)\in\ff\langle u_1\rangle$ and a set of proper parametric equations   
$$x_1=\frac{F_1(v)}{G_1(v)},\ldots,x_n=\frac{F_n(v)}{G_n(v)}$$ such that  $\mathbb I\Big(\frac{P_1(u_1)}{Q_1(u_1)},\ldots, \frac{P_n(u_1)}{Q_n(u_1)}\Big)=\mathbb I\Big(\frac{F_1(v)}{G_1(v)},\ldots, \frac{F_n(v)}{G_n(v)}\Big)$, that is, the new parametric equations have the same implicit variety as  (\ref{eq-defpdrpes}). 

\vskip3pt
The following theorem solves the proper re-parameterization problem for $\Delta$-curves.

\begin{theorem} Let (\ref{eq-defpdrpes}) be a set of PDRPEs such that each $\frac{P_i(u_1)}{Q_i(u_1)}\notin \ff.$ 
  If  (\ref{eq-defpdrpes}) is not proper, then a necessary and sufficient condition for (\ref{eq-defpdrpes}) to have a  proper re-parameterization is that $l_1=1$ in (\ref{eq-charset5.3}).
In the affirmative case, based on  {Algorithm Partial-Diff-L\"{u}roth}, we can find a new parameter $v=f(u_1)/g(u_1)$ with $f,g\in\ff\{u_1\}$ and a proper re-parameterization of (\ref{eq-defpdrpes}) 
\begin{equation}
x_1=\frac{F_1(v)}{G_1(v)},\,\ldots,\,x_n=\frac{F_n(v)}{G_n(v)}
\end{equation}
in terms of $v$.
 \end{theorem}
\proof  Suppose  (\ref{eq-defpdrpes}) is not proper and it has a proper re-parameterization 
\begin{equation} \label{eq-proper-repara} x_1=\frac{M_1(v)}{N_1(v)},\ldots, x_{n}=\frac{M_n(v)}{N_n(v)}
\end{equation} where $M_i(v), N_i(v)\in\ff\{v\}$
 are coprime for each $i$.
Then regarded as $\Delta$-ideals  in   $\ff\{x_1,\ldots,x_n\}$,  $\mathbb I\Big(\frac{P_1(u_1)}{Q_1(u_1)},\ldots, \frac{P_n(u_1)}{Q_n(u_1)}\Big)=\mathbb I\Big(\frac{M_1(v)}{N_1(v)},\ldots,  \frac{M_n(v)}{N_n(v)}\Big)$.
Since (\ref{eq-proper-repara}) is proper, 
$\ff\langle\frac{M_1(v)}{N_1(v)},\ldots, \frac{M_n(v)}{N_n(v)}\rangle=\ff\langle v\rangle$.
So there exist coprime polynomials $T_1, T_2\in\ff\{x_1,\ldots,x_n\}$ such that 
$$v=\frac{T_1\Big(\frac{M_1(v)}{N_1(v)},\ldots, \frac{M_n(v)}{N_n(v)}\Big)}{T_2\Big(\frac{M_1(v)}{N_1(v)},\ldots, \frac{M_n(v)}{N_n(v)}\Big)}.$$
 Set $$R(u_1)=\frac{T_1\Big( \frac{P_1(u_1)}{Q_1(u_1)},\ldots, \frac{P_n(u_1)}{Q_n(u_1)}\Big)}{T_2\Big(\frac{P_1(u_1)}{Q_1(u_1)},\ldots, \frac{P_n(u_1)}{Q_n(u_1)}\Big)}\in\mathcal G:=\ff\Big\langle \frac{P_1(u_1)}{Q_1(u_1)},\ldots, \frac{P_n(u_1)}{Q_n(u_1)}\Big\rangle.$$
We now show that  $R(u_1)$ is a L\"{u}roth generator of $\ff\subsetneq\mathcal G\subseteq\ff\langle u_1\rangle.$
For each $i=1,\ldots,n,$  let $F_i=\Big[x_i\cdot N_i\big(\frac{T_1(x_1,\ldots,x_n)}{T_2(x_1,\ldots,x_n)}\big)-M_i\big(\frac{T_1(x_1,\ldots,x_n)}{T_2(x_1,\ldots,x_n)}\big)\Big]\cdot T_2^\ell\in \ff\{x_1,\ldots,x_n\}$ with some suitable $\ell\in\mathbb N$. 
Since $F_i\Big(\frac{M_1(v)}{N_1(v)},\ldots,  \frac{M_n(v)}{N_n(v)}\Big)=0$, $F_i\in \mathbb I\Big(\frac{M_1(v)}{N_1(v)},\ldots,  \frac{M_n(v)}{N_n(v)}\Big)=\mathbb I\Big(\frac{P_1(u_1)}{Q_1(u_1)},\ldots, \frac{P_n(u_1)}{Q_n(u_1)}\Big)$.
So $F_i\Big(\frac{P_1(u_1)}{Q_1(u_1)},\ldots, \frac{P_n(u_1)}{Q_n(u_1)}\Big)=0$.
Thus, $\frac{P_i(u_1)}{Q_i(u_1)}= \frac{M_i(R(u_1))}{N_i(R(u_1))}\in\ff\langle R(u_1)\rangle$
and $\mathcal G=\ff\langle R(u_1)\rangle$ follows.
That is, $\mathcal G$ has a  L\"{u}roth generator.
By Theorem \ref{th-lurothcriterion}, $l_1=1.$

\vskip5pt
Conversely, suppose $l_1=1$. 
We proceed to construct a new parameter $v$ and a proper re-parameterization for  (\ref{eq-defpdrpes}).
By Theorem \ref{th-lurothcriterion}, $l_1=1$ implies that a L\"{u}roth generator exists for $\ff\langle \mathcal P(u_1)\rangle\subset\ff\langle u_1\rangle$.
Let $$B_{1,1}=\sum_{i=1}^qa_{i}(x_1,\ldots,x_n)M_i(u_1)$$ where the $M_i$ are distinct differential monomials in $u_1$. 
Following {Algorithm Partial-Diff-L\"{u}roth}, there exists some $i_0$ s.t. 
\begin{equation} \label{eq-v}
v=\frac{a_{i_0}(\frac{P_1(u_1)}{Q_1(u_1)},\ldots, \frac{P_n(u_1)}{Q_n(u_1)})}{a_1(\frac{P_1(u_1)}{Q_1(u_1)},\ldots, \frac{P_n(u_1)}{Q_n(u_1)})}\in\ff\langle u_1\rangle\backslash\ff,
\end{equation}  and  we have $\ff\langle \mathcal P(u_1)\rangle=\ff\langle v\rangle$.
Write $v$ in the reduced form as $v=\frac{f(u_1)}{g(u_1)}$ with $f, g\in\ff\{u_1\}$ and $\gcd(f,g)=1$.

\vskip3pt
Let $\mathcal C:=Q_1(u_1)x_1-P_1(u_1), \ldots,Q_n(u_1)x_n-P_n(u_1), g(u_1)y-f(u_1)$.
Clearly, $\mathcal C$ is a $\Delta$-irreducible and coherent autoreduced set under the elimination ranking $u_1\prec x_1\prec\cdots\prec x_n\prec y$ and $\sat(\C)$ is a prime $\Delta$-ideal in $\ff\{u_1,x_1,\ldots,x_n,y\}$ with a generic point $\zeta=\big(u_1,\frac{P_1(u_1)}{Q_1(u_1)},\ldots, \frac{P_n(u_1)}{Q_n(u_1)},v\big)$.
By performing Wu-Ritt's zero decomposition algorithm on $\mathcal C$ and $DS=\{Q_1,\ldots, Q_n,g\}$ w.r.t. an  elimination ranking $\mathscr R: y\prec x_1\prec\cdots\prec x_n\prec u_1$, we shall get irreducible and coherent autoreduced sets $\mathcal D_1,\ldots,\mathcal D_\ell$ such that 
$$\mathbb V(\C/DS)=\mathbb V(\sat(\C)/DS)=\bigcup_{i=1}^\ell\mathbb V(\sat(\mathcal D_i)/DS).$$
It is clear that there exists a unique $k$ such that $\mathcal D_k$ vanishes at $\zeta$,
 and  $\mathcal D_k$ is a characteristic set of $\sat(\C)$ under $\mathscr R$.
Since $\sat(\C)$ is of dimension 1 and each $\frac{P_i(u_1)}{Q_i(u_1)}$ is a $\Delta$-rational function of $v$, 
$\mathcal D_k$ must have the following form:
\[ \mathcal D_k:= G_1(y)x_1-F_1(y), \ldots, G_n(y)x_n-F_n(y), C_{1,1}(y,u_1),\ldots,C_{1,p}(y,u_1).
\]

We claim that
 \begin{equation}  \label{eq-repara}
x_1=\frac{F_1(v)}{G_1(v)}, \ldots, x_n=\frac{F_n(v)}{G_n(v)}
\end{equation}
is a proper re-parameterization of (\ref{eq-defpdrpes}).
Since $\sat(\C)=\sat(\mathcal D_k)$, by intersecting with $\ff\{x_1,\ldots,x_n\}$,
we get $\mathbb I\Big(\frac{P_1(u_1)}{Q_1(u_1)},\ldots, \frac{P_n(u_1)}{Q_n(u_1)}\Big)=\mathbb I\Big(\frac{F_1(v)}{G_1(v)},\ldots, \frac{F_n(v)}{G_n(v)}\Big)$.
And by (\ref{eq-v}), $a_1(x_1,\ldots,x_n)y-a_{i_0}(x_1,\ldots,x_n)\in\sat(\C)\cap\ff\{x_1,\ldots,x_n,y\}=\mathbb I(\frac{F_1(v)}{G_1(v)},\ldots, \frac{F_n(v)}{G_n(v)},v)$.
So $v\in\ff\langle \frac{F_1(v)}{G_1(v)}, \ldots, \frac{F_n(v)}{G_n(v)}\rangle$ and $\ff\langle \frac{F_1(v)}{G_1(v)}, \ldots, \frac{F_n(v)}{G_n(v)}\rangle=\ff\langle v\rangle$ follows.
Thus, (\ref{eq-repara}) is a proper re-parametrization of (\ref{eq-defpdrpes}).
\qed

\begin{example}
Let us consider the following PDRPEs
\begin{equation} \label{eq-exrep}
x_1=u^2,\,\,x_2=\delta_1(u)\delta_2(u).
\end{equation}
Let $PS:=\{x_1-u^2,\,\,x_2-\delta_1(u)\delta_2(u)\}$.
Under the elimination ranking $\mathscr R$ with $x_1\prec x_2\prec u$  and $\mathscr R$ restricted to each $\Theta(x_i)$ or $\Theta(u)$ being the canonical ranking $\mathscr R_0$, by Wu-Ritt's zero decomposition theorem, we have 
$$\mathbb V(PS)=\mathbb V\big(4x_1x_2-\delta_1(x_1)\delta_2(x_1), u^2-x_1\big/ux_1\big)\bigcup\mathbb V(x_1,x_2, u).$$
Only $\mathcal A:=4x_1x_2-\delta_1(x_1)\delta_2(x_1), u^2-x_1$ vanishes for $x_1=u^2,\,\,x_2=\delta_1(u)\delta_2(u).$
So $T=\{ u^2-x_1\}$. 
Thus, we can take a new  parameter $v=-x_1|_{x_1=u^2}=-u^2$.
Eliminate $u$ from $v+u^2, x_1-u^2,\,\,x_2-\delta_1(u)\delta_2(u)$,
we get a new set of  PDRPEs
\begin{equation}
x_1=-v,\,\, x_2=-\frac{\delta_1(v)\delta_2(v)}{4v},
\end{equation}
which is a proper repremeterization of (\ref{eq-exrep}).
\end{example}

 \section{Acknowledgements}
This work was supported by NSFC grants (12122118, 11971029 and 11688101),  the fund of Youth Innovation Promotion Association of CAS, and CAS Project for Young Scientists in Basic Research with Grant No. YSBR-008.

\end{document}